\documentclass[11pt]{article}%
\usepackage[latin1]{inputenc}
\usepackage{amsfonts,amssymb,latexsym}
\usepackage{amsmath,amsthm}
\usepackage[Sonny]{fncychap}
\usepackage{enumerate}

\usepackage[dvips]{color,graphicx}
\usepackage{pst-all}
\usepackage{pstcol}
\usepackage{amsbsy,amssymb}
\usepackage{amsmath}
\usepackage{amsfonts}
\usepackage{amssymb}
\usepackage{graphicx}%
\newtheorem{theorem}{Theorem}[section]
\newtheorem{definition}[theorem]{Definition}
\newtheorem{proposition}[theorem]{Proposition}

\newtheorem{lemma}[theorem]{Lemma}
\newtheorem{remark}{Remark}[section]

\numberwithin{equation}{section}
\def\fin { \vskip 0pt \hfill \hbox{\vrule height 5pt width 5pt depth 0pt} \vskip 12pt}
\begin{document}

\title{On the Schr\"{o}dinger equation with singular potentials}
\author{{{Jaime Angulo Pava}}\\{\small IME-USP, Rua do Matao 1010, Cidade Universitaria}\\{\small {CEP 05508-090, Sao Paulo, SP, Brazil.}}\\{\small \texttt{E-mail:angulo@ime.usp.br}}\vspace{.5cm} \\{{Lucas C. F. Ferreira}}\\{\small Universidade Estadual de Campinas, IMECC - Departamento de
Matem\'{a}tica,} \\{\small {Rua S\'{e}rgio Buarque de Holanda, 651, CEP 13083-859, Campinas-SP,
Brazil.}}\\{\small \texttt{E-mail:lcff@ime.unicamp.br}}}
\date{}
\maketitle

\begin{abstract}
We study the Cauchy problem for the non-linear Schr\"odinger equation with singular
potentials. For point-mass potential and nonperiodic case, we prove existence
and asymptotic stability of global solutions in weak-$L^{p}$ spaces. Specific interest is give to the point-like $\delta$ and $\delta'$ impurity and for two $\delta$-interactions in one dimension. We also consider the
periodic case which is analyzed in a functional space based on Fourier transform and
local-in-time well-posedness is proved.

\end{abstract}


{\small \vspace{0.2cm}\quad\textbf{Keywords} NLS-Dirac equation, Singular
potential, Existence, Asymptotic behavior}

\section{Introduction}

\bigskip We are interested  in this paper in the the Cauchy problem for the following Schr\"odinger model
\begin{equation}
\left\{
\begin{aligned}
i\partial_{t}u+\Delta u+\mu(x)u  &  = F(u),\quad x\in \mathbb R^n, \;\;t\in \mathbb R\ \label{SCH0}\\
u(x,0)  &  =u_{0}(x),
\end{aligned}\right.
\end{equation}
in  weak-$L^p$ spaces (Marcinkiewicz spaces) and in a space based on Fourier transform. In the weak-$L^p$ spaces we consider the  case $n=1$ and $\mu(x)=\sigma \delta$, $\mu(x)=\sigma(\delta(\cdot-a)+\delta(\cdot+a))$ (two Dirac's $\delta$ potentials place at the points $\pm a \in \mathbb R$) or   $\mu(x)=\sigma \delta'$ where $\delta$ represents the delta function in the origin and $\sigma\in \mathbb R$,  $F(u)=  \lambda\left\vert u\right\vert ^{\rho
-1}u$, where $\lambda=\pm1$ and $\rho>1$. In the space based on Fourier transform we consider $n$ arbitrary and  $\mu(x)$ being a bounded continuous function with a Fourier transform being a finite Radon measure and $F(u)= \lambda u^{\rho}$, where $\lambda=\pm1$ and $\rho\in \mathbb N$. The case $F(u)=  \lambda\left\vert u\right\vert ^{\rho
-1}u$ is also commented.

The non-linear Schr\"odinger model  (\ref{SCH0})  in the case $\mu(x)=\sigma\delta(x)$  (called the non-linear Schr\"odinger  equation with a $\delta$-type impurity, the NLS-$\delta$ equation henceforth) arise in different areas of quantum field theory and are essential for understanding a number of phenomena in condensed matter physics. At the experimental side, the recent interest in point-like impurities (defects) is triggered by the great progress in building nanoscale devises. More exactly, the NLS-$\delta$ model with a impurity at the origin in  the repulsive ($\sigma <0$) case and in the attractive ($\sigma >0$) is described by the following boundary problem (see Caudrelier\&Mintchev\&Ragoucy \cite{CMR})
\begin{equation}
\left \{
\begin{aligned}
i\partial_{t}u(x,t)+ u_{xx}(x,t) &  =  \lambda \left\vert u(x,t)\right\vert ^{\rho
-1}u(x,t),\quad x\neq0 \\
\lim_{x\to 0^+}[u(x,t)-u(-x, t)]&=0, \\
\lim_{x\to 0^+}[\partial_x u(x,t)-\partial_x u(-x, t)]&=\sigma u(0,t) \label{SCH02}\\
\lim_{x\to \pm \infty} u(x,t)=0,
\end{aligned}\right.
\end{equation}
hence $u(x,t)$ must be solution of the non-linear Schr\"odinger equation on $\mathbb R^{-}$ and $\mathbb R^{+}$, continuous at $x=0$ and satisfy a ``jump condition'' at the origin and it also vanishes at infinity.

The equations in (\ref{SCH02}) are a particular case of a more general model considering that the impurity is localized at $x=0$; in fact the equation of motion
$$
i\partial_{t}u(x,t)+ u_{xx}(x,t) =  \lambda \left\vert u(x,t)\right\vert ^{\rho
-1}u(x,t),\quad x\neq0,
$$
with the impurity boundary conditions
\begin{equation}\label{bc}
\left(\begin{array}{c}u(0+,t) \\ \partial_xu(0+,t)\end{array}\right)=\alpha
\left(\begin{array}{cc} a & b\ \\c & d\end{array}\right)\left(\begin{array}{c} u(0-,t) \\ \partial_xu(0-,t) \end{array}\right)
\end{equation}
with
\begin{equation}\label{para}
\{a,b,c,d\in \mathbb R, \alpha\in \mathbb C:  ad-bc=1,  |\alpha |=1\}.
\end{equation}

The equation (\ref{bc}) captures the interaction of the ``field'' $u$ with the impurity \cite{CMR2}. The parameters
in (\ref{para}) label  the self-adjoint extensions of the (closable) symmetric operator $H_0=-\frac{d^2}{dx^2}$ defined on the space $C_0^{\infty}(\mathbb R-\{0\})$ of smooth functions with compact support separated from the origin $x=0$. In fact, by von Neumann-Krein's theory of self-adjoint extensions for symmetric operators on Hilbert spaces, it is not difficult to show that there is a 4-parameter family of self-adjoint operators which describes all one point interactions in one-dimension of the second derivative operator $H_0$. Such a family
can be equivalently described through the family of boundary conditions at the origin
\begin{equation}\label{bc1}
\left(\begin{array}{c}\psi(0+) \\ \psi'(0+)\end{array}\right)=\alpha
\left(\begin{array}{cc} a & b\ \\c & d\end{array}\right)\left(\begin{array}{c} \psi(0-) \\ \psi'(0-)\end{array}\right)
\end{equation}
with $a,b,c,d$ and $\alpha$ satisfying the conditions in (\ref{para}) (see Theorem 3.2.3 in \cite{ak}).

Here we are interested in two specific choices of the parameters in (\ref{para}), which are relevant in physics applications (see \cite{CMR2}-\cite{CMR}). The first choice $\alpha=a=d=1$, $b=0$, $c=\sigma\neq 0$ corresponds to the case of a pure Dirac $\delta$ interaction of strength $\sigma$  (see  Theorem \ref{self}  below). The second one $\alpha=a=d=1$, $c=0$, $b=\beta\neq 0$ corresponds to the case of the so-called $\delta'$ interaction of strength $\beta$ (see  Theorem \ref{selfd}  below).

In section 2 below for convenience of the reader we present  a precise formulation for the  point interaction determined by the formal linear differential operator
\begin{equation}\label{deltaop0}
-\Delta_{\sigma}=-\frac{d^2}{dx^2}+\sigma \delta,
\end{equation}
which will be match with the singular boundary condition in (\ref{SCH02}) at every time $t$.

Existence and uniqueness of local and global-in-time solutions of problem (\ref{SCH0}) with $\mu(x)=0$ and $F(u)=  \lambda\left\vert u\right\vert ^{\rho-1}u$ have been much studied in the framework of the Sobolev spaces $H^s(\mathbb R^n)$, $s\geqq 0$, i.e, the solutions and their derivatives have finite energy (see Cazenave's book \cite{C1} and the reference therein). In the case of $\delta$-interaction, namely, $\mu(x)=\sigma \delta$ the existence of global solution in $H^1(\mathbb R)$ and $L^2(\mathbb R)$ has been addressed in Adami\&Noja \cite{Adami} (we can also to apply Theorem 3.7.1 in  \cite{C1} for obtaining a local-in-time well-posedness theory in $H^1(\mathbb R)$).

The first study of infinite $L^2$-norm solutions for $\mu(x)=0$ and $F(u)=\lambda\left\vert u\right\vert ^{\rho-1}u$ was addressed by Cazenave\&Weissler in \cite{Cazenave1} where they consider the space
$$
X_\rho=\{u:\mathbb{R}\rightarrow L^{\rho+1}(\mathbb R^n) \text{ \ Bochner meas.}: \sup_{-\infty<t<\infty}|t|^{\vartheta }\|u(t)\|_{L^{\rho+1}}<\infty\},
$$
where $\vartheta =\frac{1}{\rho-1}-\frac{n}{2(\rho+1)}$ and  $\|\cdot\|_{L^{\rho+1}}$ denotes the usual ${L^{\rho+1}}$ norm. Under a suitable smallness condition on the initial data, they prove the existence of global solution in $X_\rho$, for $\rho_{0}(n)<\rho<\gamma(n)$ where $\rho_{0}(n)=\frac{n+2+\sqrt{n^{2}+12n+4}}{2n}>1$ is the positive root of the equation $n\rho^{2}-(n+2)\rho-2=0$ and $\gamma(n)=\infty$ if $n=1,2$ and $\gamma(n)=$ $\frac{n+2}{n-2}$ in otherwise.

Later on, in Cazenave\&Vega\&Vilela \cite{Cazenave2} the Cauchy problem was studied in the framework of weak-$L^p$ spaces. Using a Strichartz-type inequality, the authors obtained existence of solutions in the class $L^{(p,\infty)}(\mathbb R^{n+1})\equiv L_t^{(p,\infty)}(L_x^{(p,\infty)})$, where $(x,t)\in \mathbb R^{n}\times \mathbb R$ and $p=\frac{(\rho-1)(n+2)}{2\rho}$, for $\rho$ in the range
$$
\rho_0<\frac{4(n+1)}{n(n+2)}<\rho-1<\frac{4(n+1)}{n^2}<\frac{4}{n-2}.
$$
More recently, in Braz e Silva\&Ferreira\&Villamizar-Roa \cite{BFV} the Cauchy problem was studied in the Marcinkiewicz space $L^{(\rho+1,\infty)}$. Using bounds for the Schr\"odinger linear group in the context of Lorentz spaces, the authors showed existence and uniqueness of local-in-time solutions in the class
$$
\{u:\mathbb{R}\rightarrow L^{(\rho+1,\infty)} \text{ \ Bochner
meas.}: \sup_{-T<t<T}|t|^{\zeta}\|u(t)\|_{L^{(\rho+1,\infty)}}<\infty\},
$$
where $1<\rho<\rho_{0}(n)$ and $\frac{n(\rho-1)}{2(\rho+1)}=\zeta_{0}<\zeta<\frac{1}{\rho}$.
Since $\rho_{0}(n)<\frac{4}{n}$, the range for $\rho$ is different from the ones in Cazenave\&Weissler \cite{Cazenave1} and Cazenave {\it et al.} \cite{Cazenave2}. The existence of global solutions is showed in norms of type $\sup_{|t|>0}|t|^{\vartheta }\|u(t)\|_{L^{(\rho+1,\infty)}}$, where $\vartheta =\frac{1}{\rho-1}-\frac{n}{2(\rho+1)}$ and $\rho_0(n)<\rho<\gamma(n)$.

Our approach is based in some ideas in \cite{BFV}, so via real interpolation techniques we establish bounds for the Schr\"odinger linear group $G_\sigma(t)=e^{i(\partial_x^2+\sigma \delta)t}$ in the context of Lorentz spaces in the one-dimensional case. The cases $n=2,3$ remain open. The fundamental solution of the corresponding linear time-dependent Schr\"odinger  equation, namely
$$
iu_t=-(\Delta+\sigma \delta)u,
$$
is now well know for $ n=1,2,3$; see Albeverio {\it et al.} \cite{ABD}- for instance. However, surprisingly, a \textquotedblleft good formula\textquotedblright of the unitary group $G_\sigma(t)\phi=e^{i(\Delta+\sigma \delta)t}\phi$ depending of the free linear propagator $e^{i\Delta t}\phi$ was found explicitly only for the one-dimensional case (see Holmer {\it et al.} \cite{Holmer5}-\cite{Holmer3}). In fact, by using scattering techniques, it was established in \cite{Holmer5} the convenient formula (for the case $\sigma\geqq 0$)
\begin{equation}\label{pospro1}
G_\sigma(t) \phi(x)= e^{it \Delta} (\phi\ast \tau_\sigma) (x) \chi^0_{+} + \Big[e^{it \Delta} \phi(x) + e^{it \Delta} (\phi\ast \rho_\sigma) (-x) \Big ]\chi^0_{-}
\end{equation}
where
$$
\rho_\sigma(x)=-\frac{\sigma}{2} e^{\frac{\sigma}{2} x}\chi^0_{-},\;\; \tau_\sigma(x)=\delta (x)+ \rho_\sigma(x),
$$
with $\chi^0_{+}$ the characteristic function of $[0,+\infty)$ and $\chi^0_{-}$ the characteristic function of $(-\infty, 0]$. For the case $\sigma<0$ see \cite{Holmer3} and Theorem \ref{expli} below. Here we show in the Appendix how to obtain the formula (\ref{pospro1}) based in the fundamental solution found in Albeverio {\it et al.} \cite{ABD}, which does not use scattering ideas. Nice formulas as that in (\ref{pospro1}) are not known for the cases $n=2,3$. Now, in the case $\sigma\geqq 0$, one can show from (\ref{pospro1}) the dispersive estimate in Lorentz spaces (see Lemma \ref{grupint1} below)
\begin{equation}\label{grupint02}
\left\Vert G_\sigma(t)f\right\Vert _{(p^{\prime},d)}\leq C|t|^{-\frac
{1}{2}(\frac{2}{p}-1)}\left\Vert f\right\Vert _{(p,d)},
\end{equation}
for $1\leqq d\leqq \infty$, $p'\in (2,\infty)$, $p\in (1,2)$ and $p'$ satisfying $\frac{1}{p}+\frac{1}{p'}=1$, where $C>0$ is independent of $f$ and $t\neq 0$. Then, under a suitable smallness condition on the initial data $u_0$, the existence of global solutions for (\ref{SCH0}) is proved in the space (see Theorem \ref{GlobalTheo})
$$
\mathcal{L}_{\vartheta }^{\infty}=\{u:\mathbb{R}\rightarrow L^{(\rho+1,\infty)} \text{ \ Bochner
meas.}: \sup_{-\infty<t<\infty}|t|^{\vartheta }\Vert u\Vert_{(\rho+1,\infty)}<\infty\},
$$
for $\sigma\geq0$, where $\vartheta =\frac{1}{\rho-1}-\frac{1}{2(\rho+1)}$ and $\rho_{0}=\frac{3+\sqrt{17}}{2}>1$ is  the positive root of the equation $\rho^{2}-3\rho-2=0$. We also analyze the asymptotic stability of the global solutions (see Theorem \ref{TeoAssin}). For $\sigma <0$ our approach in general is not applicable because in this case the operator $-\Delta_\sigma$ has a non-trivial negative point spectrum. But, in this case it is possible to show the existence of a invariant manifold of periodic orbits in Lorentz spaces (see Section 6).

With regard to the more two singular cases: two Dirac's $\delta$ potentials placed at the points $\pm a \in \mathbb R$, $\mu(x)=\alpha(\delta(\cdot-a)+\delta(\cdot+a))$, and $\mu(x)=\beta \delta'$, i.e., the derivative of a $\delta$, a similar analysis to that above for the case of a $\delta$-potential  can be established. In these cases, it is not known an explicit expression for the associated time propagator as that in (\ref{pospro1}) for the case of $\mu(x)=\sigma\delta$. However, by using a formula for the integral kernel of the time propagator associated (see \cite{KS} and \cite{AGHH}), we obtain an estimate similar to (\ref{grupint02}).

For the case $n\geqq 4$ we do not have Schr\"odinger operators with point interactions. In fact, the Schr\"odinger operators with point interactions, namely, perturbations of the Laplace operators by ``measures'' supported on a discrete set (supported at zero for simplicity, namely, by the Dirac delta measure $\delta$ centered at zero) are usually defined by means of von Neumann\&Krein theory of self-adjoint extensions of symmetric operators, and so as one of a whole family of self-adjoint (in $L^2(\mathbb R^n)$) extensions of an operator $A$, $D(A)=C^\infty_0(\mathbb R^n-\{0\})$, $Au=-\Delta u$, $u\in D(A)$. In the case $n\geqq 4$  it is well known that the theory trivializes where there is only one self-adjoint extension of $A$ (see Albeverio {\it et al.} \cite{AGHH}).

Next, let $\mathcal M$ be the set of finite Radon measure endowed with the norm of total variation, that is, $\|\omega\|_{\mathcal M}=|\omega|(\mathbb R^n)$ for $\omega\in \mathcal M$, $n\geq1$. Then, by considering  $\mu(x)$ in (\ref{SCH0}) being a bounded continuous function with a Fourier transform such that $\widehat{\mu}\in \mathcal M$, we show a local-in-time well-posedness result in the Banach space
\begin{equation}
\mathcal{I}=[\mathcal{M}(\mathbb{R}^{n})]^{\vee}=\{f\in\mathcal{S}^{\prime
}(\mathbb{R}^{n}):\widehat{f}\in\mathcal{M}(\mathbb{R}^{n})\}
\end{equation}
whose norm is given by $\|f\|_\mathcal I=\|\widehat{f}\|_{\mathcal M}$. We also obtain a similar result in the periodic case (see Section 7).

\section{The one-center  $\delta$-interaction in one dimension }

In this subsection for convenience of the reader  we establish initially a precise formulation for the  point interaction determined by the formal linear differential operator
\begin{equation}\label{deltaop}
-\Delta_{\sigma}=-\frac{d^2}{dx^2}+\sigma \delta,
\end{equation}
defined on functions on the real line. The parameter $\sigma$ represents the coupling constant or strength attached to the point source located at $x=0$. We note that there are many approaches for studying the operator in (\ref{deltaop}), for instance, by the use of quadratic forms or by the self-adjoint extensions of symmetric operators. We also note that the quantum mechanics model in (\ref{deltaop}) has been studied into a more general framework when it is associated with the Kronig-Penney model in solid state physics (see Chapter III.2 in Albeverio {\it et al.} \cite{AGHH}) or when it is associated to singular rank one perturbations (Albeverio {\it et al.} \cite{ak}).

By following \cite{ak}, we consider the operator $A=-\frac{d^2}{dx^2}$ with domain $D(A)= H^2(\mathbb R)$  and the (closeable)  symmetric restriction $A^0\equiv A|_{D(A^0)}$ with dense domain $D(A^0)=\{\psi\in D(A) : (\delta, \psi)\equiv\psi(0)=0\}$. Then we obtain that the {\it deficiency subspaces} of $A^0$,
 \begin{equation}\label{sa1}
\mathcal D_+=\text{Ker}({A^0}^*-i),\quad{\rm{and}}\quad\mathcal D_-=\text{Ker}({A^0}^*+i),
\end{equation}
have dimension ({\it deficiency indexes}) equal to $1$. It is no difficult to see that these subspaces are generated, respectively, by  $g_{+i}\equiv (A-i)^{-1}\delta$ and  $g_{-i}\equiv (A+i)^{-1}\delta$,  called {\it deficiency elements}  and given explicitly by (see \cite{ak}),
 \begin{equation}\label{def}
g_{\pm i}(x)=\frac{i}{2\sqrt{\pm i}}e^{i\sqrt{\pm i} |x|},\qquad Im \sqrt{\pm i} >0.
\end{equation}
We note that the Fourier transform of $g_{\pm i}$ are given by $\widehat{g_{\pm i}}(\xi)=\frac{1}{\xi^2\mp i}$

Next we present explicitly  all the self-adjoint extensions of the symmetric operator $A^0$, which will be parameterized by the strength $\sigma$. By normalizing the deficiency elements $\widetilde{g}_{\pm i}=\frac{g_{\pm i}}{\|g_{\pm i\|}}$ and for convenience of notation we will continue to use $g_{\pm i}$, we have from the von Neumann's theory of self-adjoint extensions for symmetric operators (see \cite{RS})  that all the closed symmetric extensions of $A^0$ are self-adjoint and coincides with the restriction of the operator ${A^0}^*$. Moreover, for $\theta\in[0,2\pi)$ the self-adjoint extension $A^0(\theta)$ of $A^0$ is defined as follows:
\begin{equation}\label{sa4}
\left\{
\begin{aligned}
&D(A^0(\theta))=\{\psi+\lambda g_i+\lambda e^{i\theta}g_{-i} :\psi\in D(A^0), \lambda\in \mathbb C\},\\
&A^0(\theta)(\psi+\lambda g_i+\lambda e^{i\theta}g_{-i})={A^0}^*(\psi+\lambda g_i+\lambda e^{i\theta}g_{-i})=A^0\psi+i\lambda g_i-i\lambda e^{i\theta}g_{-i}.
\end{aligned}\right.
\end{equation}
Thus from \eqref{sa4} and  \eqref{def} we obtain that for $\zeta\in D(A^0(\theta))$,  in the form $\zeta=\psi+\lambda g_i+\lambda e^{i\theta}g_{-i}$, we have the basic expression
\begin{equation}\label{sa5}
\zeta'(0+)-\zeta'(0-)=-\lambda(1+e^{i\theta}).
\end{equation}
Next we find $\sigma$ such that $\sigma\zeta(0)=-\lambda(1+e^{i\theta})$. Indeed, $\sigma$ is given by the formula
\begin{equation}\label{sa6}
\sigma(\theta)=\frac{-2\cos(\theta/2)}{cos(\frac{\theta}{2}-\frac{\pi}{4})}.
\end{equation}

So, from now on we parameterize all self-adjoint extensions of $A^0$ with the help of $\sigma$. Thus we get:

\begin{theorem}\label{self} All self-adjoint extensions of $A^0$ are given for $-\infty<\sigma\leqq +\infty$ by
\begin{equation}\label{sa8}
\left\{
\begin{aligned}
-\Delta_{\sigma}&=-\frac{d^2}{dx^2}\\
D(-\Delta_{\sigma})&=\{\zeta\in H^1(\mathbb R)\cap H^2(\mathbb R-\{0\}):  \zeta'(0+)-\zeta'(0-)=\sigma \zeta(0)\}.
\end{aligned}\right.
\end{equation}
The special case $\sigma=0$ just leads to the operator $-\Delta$ in $L^2(\mathbb R)$,
\begin{equation}\label{sa9}
-\Delta=-\frac{d^2}{dx^2},\qquad D(-\Delta)= H^2(\mathbb R),
\end{equation}
whereas the case $\sigma=+\infty$ yields a Dirichlet boundary condition at zero,
\begin{equation}\label{sa9a}
D(-\Delta_{+\infty})=\{\zeta\in H^1(\mathbb R)\cap H^2(\mathbb R-\{0\}):  \zeta(0)=0\}.
\end{equation}

\end{theorem}

\textbf{Proof.} By the arguments sketched above we obtain easily that $A^0(\theta)\subset -\Delta_{\sigma}$,
with $\sigma=\sigma(\theta)$ given in \eqref{sa6}. But $-\Delta_{\sigma}$ is  symmetric in the corresponding domain $D(-\Delta_{\sigma})$ for all $-\infty<\sigma\leqq +\infty$, which implies the relation $A^0(\theta)\subset -\Delta_{\sigma}\subset (-\Delta_{\sigma})^* \subset A^0(\theta)$.
It completes the proof of the Theorem.
\fin

Next, we recall the basic spectral properties of $-\Delta_\sigma$ which will be relevant for  our  results  (see \cite{AGHH}).

\begin{theorem}\label{resol5b} Let $-\infty<\sigma\leqq +\infty$. Then the essential spectrum of $-\Delta_\sigma$ is the nonnegative real axis, $\Sigma_{ess}(-\Delta_\sigma)=[0,+\infty)$.

 If $-\infty<\sigma<0$, $-\Delta_\sigma$ has exactly one negative, simple eigenvalue, i.e., its discrete spectrum $\Sigma_{dis}(-\Delta_\sigma)$ is $\Sigma_{dis}(-\Delta_\sigma)=\{{-\sigma^2/4}\}$, with a strictly (normalized) eigenfunction
 $$
 \Psi_\sigma(x)=\sqrt{\frac{-\sigma}{2}}e^{\frac{\sigma}{2}|x|}.
 $$

 If $\sigma\geqq 0$ or $\sigma=+\infty$, $-\Delta_\sigma$ has not discrete spectrum, $\Sigma_{dis}(-\Delta_\sigma)=\emptyset$.

\end{theorem}

\section{Two symmetric  $\delta$-interaction in one dimension }

The  one-dimensional Schr\"odinger operator with two symmetric delta interactions of strength $\alpha$ and placed at the point $\pm a$ is given formally by the  linear differential operator
\begin{equation}\label{deltaop1}
-\Delta_{\alpha}=-\frac{d^2}{dx^2}+\alpha( \delta(\cdot-a)+\delta(\cdot+a)),
\end{equation}
defined on functions on the real line.  By using the same notations as in last section, the symmetric operator $A^1=A|_{D(A^1)}$ with dense domain
$$
D(A^1)=\{\psi\in H^2(\mathbb R): \psi(\pm a)=0\},
$$
has deficiency indices (2, 2), and so from the Von Neumann-Krein theory we have that all self-adjoint extensions of $A^1$ are given by a four-parameter family of self-adjoint operators. Here we restrict to the case of so-called separated boundary conditions at each point $\pm a$. More specifically, we have the following theorem (see \cite{ABD}).

\begin{theorem}\label{self3} There is a family of self-adjoint extensions of $A^1$  given for $-\infty<\alpha\leqq +\infty$ by
\begin{equation}\label{sa8}
\left\{
\begin{aligned}
-\Delta_{\alpha}&=-\frac{d^2}{dx^2}\\
D(-\Delta_{\alpha})&=\{\zeta\in H^1(\mathbb R)\cap H^2(\mathbb R-\{\pm a\}):  \zeta'(\pm a+)-\zeta'(\pm a-)=\alpha \zeta(\pm a)\},
\end{aligned}\right.
\end{equation}
 The special case $\alpha=0$ just leads to the operator $-\Delta$ in $L^2(\mathbb R)$,
\begin{equation}\label{sa9}
-\Delta=-\frac{d^2}{dx^2},\qquad D(-\Delta)= H^2(\mathbb R),
\end{equation}
whereas the case $\alpha=+\infty$ yields a Dirichlet boundary condition at  the point $\pm a$,
\begin{equation}\label{sa9a}
D(-\Delta_{+\infty})=\{\zeta\in H^1(\mathbb R)\cap H^2(\mathbb R-\{\pm a\}):  \zeta( \pm a+)=\zeta( \pm a-)=0\}.
\end{equation}

\end{theorem}

Next, we establish the basic spectral properties of $-\Delta_\alpha$ which will be relevant for  our  results  (see \cite{AGHH}).

\begin{theorem}\label{resol2de} Let $-\infty<\alpha\leqq +\infty$. Then the essential spectrum of $-\Delta_\alpha$ is the nonnegative real axis, $\Sigma_{ess}(-\Delta_\alpha)=[0,+\infty)$.

 I) If $-\infty<\alpha<0$, then the discrete spectrum of $-\Delta_\alpha$, $\Sigma_{dis}(-\Delta_\alpha)$, consists of negative eigenvalues $\gamma$ given by the implicit equation
 $$
 (-2i \eta +\alpha)^2=\alpha^2e^{4i\eta a},\quad \eta=\sqrt{\gamma},\;\; Im \;\eta>0.
 $$
Moreover, we have that:
 \begin{enumerate}

 \item[1)] if $a\leqq -\frac{1}{\alpha}$, then $\Sigma_{dis}(-\Delta_\alpha)=\{\gamma_1(a,\alpha)\}$, where $\gamma_1(a,\alpha)$ is defined by
 $$
\gamma_1(a,\alpha)=-\frac{1}{4a^2}[W(-a\alpha e^{a\alpha})-a\alpha]^2,
$$
 where $W(\cdot)$ is the Lambert special function (or product logarithm) defined by the equation $W(x)e^{W(x)}=x$.

  \item[2)]  if $a> -\frac{1}{\alpha}$, then $\Sigma_{dis}(-\Delta_\alpha)=\{\gamma_1(a,\alpha), \gamma_2(a,\alpha)\}$, where $\gamma_2(a,\alpha)$ is defined by
 $$
\gamma_2(a,\alpha)=-\frac{1}{4a^2}[W(a\alpha e^{a\alpha})-a\alpha]^2.$$
 \end{enumerate}

II)  If $\alpha\geqq 0$ or $\alpha=+\infty$, $-\Delta_\alpha$ has not discrete spectrum, $\Sigma_{dis}(-\Delta_\alpha)=\emptyset$.

\end{theorem}

\section{The  $\delta'$-interaction in one dimension }

In this subsection for convenience of the reader  we establish  a precise formulation for the  point interaction determined by the formal linear differential operator
\begin{equation}\label{deltaop1}
-\Delta_{\beta}=-\frac{d^2}{dx^2}+\beta \delta',
\end{equation}
defined on functions on the real line. The parameter $\beta$ represents the coupling constant or strength attached to the point source located at $x=0$ and $ \delta'$ is the derivative of the $\delta$. By following Albeverio {\it et al.} \cite{AGHH}-\cite{ak}, the elements in the domain of the operator $-\Delta_{\beta}$ are characterized by suitable bilateral singular boundary conditions at the singularity (see (\ref{bc1})),  while the real true action coincides with the laplacian out the singularity. At variance with the $\delta$ interaction, whose domain is contained in $H^1(\mathbb R)\cap H^2(\mathbb R-\{0\})$ (in particular in a continuous function set, see (\ref{sa8})), the latter has a domain  contained only in $H^2(\mathbb R-\{0\})$ and so by allowing discontinuities of the elements at the position of the defect. More precisely, for $A^2=A|_{D(A^2)}$ being considered with dense domain
$$
D(A^2)=\{\psi\in H^2(\mathbb R): \psi(0)=\psi'(0)=0\}.
$$
$A^2$ has deficiency indices $(2,2)$ and hence it has a four-parameter family of self-adjoint. We are interested in the following one-parameter family of self-adjoint extensions (see \cite{AGHH}-\cite{ak}).

\begin{theorem}\label{selfd} There is  a family of self-adjoint extensions of $A^2$  given for $-\infty<\beta\leqq +\infty$ by
\begin{equation}\label{sa8d}
\left\{
\begin{aligned}
-\Delta_{\beta}&=-\frac{d^2}{dx^2}\\
D(-\Delta_{\beta})&=\{\zeta\in H^2(\mathbb R-\{0\}):  \zeta'(0+)=\zeta'(0-), \zeta(0+)-\zeta(0-)=\beta\zeta'(0-)\},
\end{aligned}\right.
\end{equation}
 The special case $\beta=0$ just leads to the operator $-\Delta$ in $L^2(\mathbb R)$,
\begin{equation}\label{sa9d}
-\Delta=-\frac{d^2}{dx^2},\qquad D(-\Delta)= H^2(\mathbb R),
\end{equation}
whereas the case $\beta=+\infty$ yields a Neumann boundary condition at zero and decouples $(-\infty, 0)$ and $(0,+\infty)$, i.e.,
\begin{equation}\label{sa9d}
D(-\Delta_{+\infty})=\{\zeta\in H^2(\mathbb R-\{0\}):  \zeta'(0+)=\zeta'(0-)=0\}.
\end{equation}
\end{theorem}

Note that the functions in the domain of $\delta'$ have a jump at the origin, and the left and right derivatives coincide. Next, we establish the basic spectral properties of $-\Delta_\beta$ which will be relevant for  our  results  (see \cite{AGHH}).

\begin{theorem}\label{resol5d} Let $-\infty<\beta\leqq +\infty$. Then the essential spectrum of $-\Delta_\beta$ is the nonnegative real axis, $\Sigma_{ess}(-\Delta_\beta)=[0,+\infty)$.

 If $-\infty<\beta<0$, $-\Delta_\beta$ has exactly one negative simple eigenvalue, i.e., its discrete spectrum $\Sigma_{dis}(-\Delta_\beta)$ is $\Sigma_{dis}(-\Delta_\beta)=\{-\frac{4}{\beta^2}\}$, with a (normalized) eigenfunction
 $$
 \Phi_\beta(x)=\Big(-\frac{2}{\beta}\Big)^{\frac12}\text{sign}(x) e^{\frac{2}{\beta}|x|}.
 $$

 If $\beta\geqq 0$ or $\beta=+\infty$, $-\Delta_\beta$ has not discrete spectrum, $\Sigma_{dis}(-\Delta_\beta)=\emptyset$.

\end{theorem}

\subsection{The linear propagator }

\subsubsection{The case $\mu(x)=\sigma \delta$.}

Next we determine the linear propagator $G_\sigma=e^{i(\Delta+\sigma\delta)t}$(unitary group) determined by the linear system associated with (\ref{SCH0}),
\begin{equation}
\left \{
\begin{aligned}
i\partial_{t}u&=-(\Delta u+\sigma\delta)u\equiv H_\sigma u,  \label{SCH1}\\
u(0)  &  =u_{0},
\end{aligned} \right.
\end{equation}
where  we are using the notation $H_\sigma=-\Delta_{-\sigma}$.

We will use the representation of the propagator in terms of the  eigenfunctions (associated to discrete eigenvalues) and generalized eigenfunctions (see Iorio \cite{Io}, Holmer {\it et al.} \cite{Holmer5} and Duch\^ene {\it et al.}  \cite{DMW}). Indeed, the family of generalized eigenfunctions $\{\psi_\lambda\}_{\lambda\in \mathbb R}$ will be such that satisfy
\begin{equation}
\left\{
\begin{aligned}\label{eigengen}
&H_\sigma \psi_\lambda=\lambda^2 \psi_\lambda,\qquad  \psi_\lambda\;\text{continuous and}\\\
& \psi'_\lambda(0+)-\psi'_\lambda(0-)=\sigma\psi_\lambda(0).
\end{aligned}\right.
\end{equation}
Hence we obtain the following family of special solutions, $e_{\pm}(x,\lambda)$ to (\ref{eigengen}), as follows
\begin{equation}\label{efun}
e_{\pm}(x,\lambda)=t_\sigma(\lambda)e^{\pm i\lambda x}\chi^0_{\pm}+(e^{\pm i\lambda x}+r_\sigma(\lambda)e^{\mp i\lambda x})\chi^0_{\mp},
\end{equation}
where  $\chi^0_{+}$ is the characteristic function of $[0,+\infty)$ and $\chi^0_{-}$ is the characteristic function of $(-\infty, 0]$. $t_\sigma$ and $r_\sigma$ are the transmission and reflection coefficients:
\begin{equation}\label{tr}
t_\sigma(\lambda)=\frac{2i\lambda}{2i\lambda-\sigma},\quad r_\sigma(\lambda)=\frac{\sigma}{2i\lambda-\sigma}.
\end{equation}
They satisfy the following two equations:
\begin{equation}\label{tr2}
| t_\sigma(\lambda)|^2+| r_\sigma(\lambda)|^2=1,\quad r_\sigma(\lambda)+1=t_\sigma(\lambda).
\end{equation}

Next, by defining the family $\{\psi_\lambda\}_{\lambda\in \mathbb R}$ as
$$
\psi_\lambda (x)=\left\{
\begin{aligned}
& e_+(x,\lambda)\;\; \;\;for  \;\; \lambda\geqq 0\\
&e_-(x, -\lambda)\;\; for \;\; \lambda < 0
\end{aligned} \right.
$$
we obtain from Theorem \ref{resol5b} the following relations (see \cite{DMW}, \cite{Io});
\begin{enumerate}\label{orthorela}
\item[1)] $\int_{\mathbb R} \Psi_\sigma(x) \overline{\psi_\lambda(x)}dx=0$,  \;\; for all $\lambda \in \mathbb R$ and $ \sigma <0$,\\
\item[2)] $ \int_{\mathbb R} \psi_\mu(x) \overline{\psi_\lambda (x)}dx=\delta (\lambda-\mu)$, \;\; for all $\mu, \lambda \in\mathbb R$,\\
\item[3)] $\Psi_\sigma(x)\Psi_\sigma(y)+  \int_{\mathbb R} \psi_\lambda(x) \overline{\psi_\lambda (y)}d\lambda=\delta (x-y)$, \;\; $ \lambda \in \mathbb R$, $\sigma <0$.
\end{enumerate}
We recall that the relation 3) above, called the {\it completeness relations}, in the case $\sigma >0$ is reads as $\int_{\mathbb R} \psi_\lambda(x) \overline{\psi_\lambda}(y)d\lambda=\delta (x-y)$ (the proof of 3) for the family $\{\psi_\lambda\}_{\lambda\in \mathbb R}$ can be showed by following the ideas in the proof of Theorem \ref{expli} below). Moreover, the family $\{\psi_\lambda\}_{\lambda\in \mathbb R}$ allows us to define the {\it generalized Fourier transform}
\begin{equation}
\mathcal F (f)(\lambda)=\int_{\mathbb R} f(x)\overline{\psi_\lambda (x)}dx,
\end{equation}
and its formal adjoint $\mathcal G (g)(x)=\int_{\mathbb R} \psi_\lambda (x) g(\lambda)d\lambda$.  Hence, from  2) we obtain immediately that $\mathcal G$ is the inverse Fourier transform, namely,
$$
f(\lambda)=f\ast \delta (\lambda)=\int_{\mathbb R} \overline{\psi_\lambda (x)}\int_{\mathbb R} f(\mu)\psi_\mu (x)d\mu dx=\mathcal F(\mathcal G g)(\lambda).
$$
Moreover, from the completeness relations 3) we obtain for every $f\in L^2(\mathbb R)$ the following  (orthogonal) expansion in eigenfunctions of $H_\sigma$,
\begin{equation}
f= \langle f, \Phi_\sigma\rangle \Phi_\sigma + \int_{\mathbb R} \mathcal F( f)(\lambda)\psi_\lambda(x)d\lambda.
\end{equation}
Thus for $u\in C(\mathbb R; L^2(\mathbb R))$ being a solution of (\ref{SCH1}), the method of separation of variables implies that
\begin{equation}\label{solu1}
u(x,t)=e^{-it H_\sigma} u_0(x)=
e^{i\frac{\sigma^2}{4}t} \langle u_0, \Phi_\sigma\rangle \Phi_\sigma (x) +\int_{\mathbb R}e^{-i\lambda^2 t}\mathcal F (u_0)(\lambda)\psi_\lambda (x)d\lambda.
\end{equation}

In the next Theorem we describe  explicitly the propagator $e^{-it H_\sigma}$ in terms of the free propagator of the Schr\"odinger equation $e^{it\Delta }$ (see Holmer {\it et al.} \cite{Holmer5}, Datchev\&Holmer \cite{Holmer3}). In the Appendix we present a different proof based in the fundamental solution associated to (\ref{SCH1}).

\begin{theorem} \label{expli} Suppose that $ \phi\in L^1(\mathbb R)$  with $supp\; \phi \subset (-\infty, 0]$. Then,
\begin{enumerate}
\item[1)] Para $\sigma \geqq 0$, we have
\begin{equation}\label{pospro}
e^{-it H_\sigma} \phi(x)= e^{it \Delta} (\phi\ast \tau_\sigma) (x) \chi^0_{+} + \Big[e^{it \Delta} \phi(x) + e^{it \Delta} (\phi\ast \rho_\sigma) (-x) \Big ]\chi^0_{-}
\end{equation}
where
$$
\rho_\sigma(x)=-\frac{\sigma}{2} e^{\frac{\sigma}{2} x}\chi^0_{-},\;\; \tau_\sigma(x)=\delta (x)+ \rho_\sigma(x).
$$

\item[ 2)] Para $\sigma < 0$, we have
\begin{equation}\label{negpro}
e^{-it H_\sigma} \phi(x)= e^{i\frac{\sigma^2}{4} t} P\phi (x)+ e^{it \Delta} (\phi\ast \tau_\sigma) (x) \chi^0_{+} + \Big[e^{it \Delta} \phi(x) + e^{it \Delta} (\phi\ast \rho_\sigma) (-x) \Big ]\chi^0_{-}
\end{equation}
where
$$
\rho_\sigma(x)=\frac{\sigma}{2} e^{\frac{\sigma}{2} x}\chi^0_{+},\;\; \tau_\sigma(x)=\delta (x)+ \rho_\sigma(x),
$$
and $P$ is the $L^2(\mathbb R)$-orthogonal projection onto the eigenfunction $\Phi_\sigma$, $P\phi=\langle \phi, \Phi_\sigma\rangle \Phi_\sigma $.
\end{enumerate}
\end{theorem}

\begin{remark} We observe the following:
\begin{itemize}
\item[1) ] The Fourier transform de $\rho_\sigma$ for every sign of $\sigma$ is given by $\widehat{\rho_\sigma}(\lambda)=r_\sigma(\lambda)$ and so $\widehat{\tau_\sigma}(\lambda)=1+r_\sigma(\lambda)=t_\sigma(\lambda)$.
\item[2) ] Formula (\ref{pospro}) and (\ref{negpro}) will allow us to estimate the operator norm of $e^{-it H_\sigma}$ using $e^{it \Delta}$.
\end{itemize}
\end{remark}

\textbf{Proof.} We only consider the case $\sigma\geqq 0$.  From (\ref{solu1}), without the term of projection, we have from the definition of the family $\{\psi_\lambda\}$ and a change of variable that
\begin{equation}\label{grupo1}
e^{-it H_\sigma} \phi(x)=\int_{\mathbb R}\phi(y)\int_{0}^\infty e^{-it\lambda^2}\Big (e_+(x,\lambda)\overline{e_+(y,\lambda)}+e_-(x,\lambda)\overline{e_-(y,\lambda)}\Big)d\lambda dy.
\end{equation}

Next, we compute first
\begin{equation}
\begin{aligned}
& \int_{\mathbb R}\phi(y)\overline{e_+(y,\lambda)}dy=\int_{-\infty}^0 \phi(y)e^{-i\lambda y}dy +\overline{r_\sigma(\lambda)}\int_{-\infty}^0 \phi(y)e^{i\lambda y}dy=\widehat{\phi}(\lambda)+r_\sigma(-\lambda)\widehat{\phi}(-\lambda),\\
&\qquad\qquad\qquad \qquad\int_{\mathbb R}\phi(y)\overline{e_-(y,\lambda)}dy=t_\sigma(-\lambda)\widehat{\phi}(-\lambda),
\end{aligned}
\end{equation}
so for $x>0$  we have from (\ref{grupo1}) and the fact $r_\sigma(-\lambda)t_\sigma(\lambda)+r_\sigma(\lambda)t_\sigma(-\lambda)=0$ that
\begin{equation}
e^{-it H_\sigma} \phi(x)=\int_{\mathbb R}e^{-it\lambda^2}t_\sigma(\lambda)\widehat{\phi}(\lambda)e^{i\lambda x}d\lambda=e^{it \Delta}(\tau_\sigma\ast\phi)(x),
\end{equation}
where $\widehat{\tau_\sigma}(\lambda)=t_\sigma(\lambda)$. Similarly, since $r_\sigma(-\lambda)r_\sigma(\lambda)+t_\sigma(-\lambda)t_\sigma(\lambda)=1$, we have for $x<0$
\begin{equation}
e^{-it H_\sigma} \phi(x)=\int_{\mathbb R}e^{-it\lambda^2}(\widehat{\phi}(\lambda)e^{i\lambda x}+r_\sigma(\lambda)\widehat{\phi}(\lambda)e^{-i\lambda x})d\lambda=e^{it \Delta}\phi(x)+e^{it \Delta}(\rho_\sigma\ast\phi)(-x),
\end{equation}
where $\widehat{\rho_\sigma}(\lambda)=r_\sigma(\lambda)$.
\fin

\subsubsection{The case $\mu(x)=\alpha (\delta(\cdot-a)+\delta(\cdot+a))$.}

Next we determine the linear propagator $M_\alpha(t)=e^{-itU_\alpha} $ (unitary group) determined by the linear system associated with (\ref{SCH0}),
\begin{equation}\label{beta3}
\left \{
\begin{aligned}
i\partial_{t}u&=-(\Delta u+\alpha (\delta(\cdot-a)+\delta(\cdot+a))u\equiv U_\alpha u, \\
u(0)  &  =u_{0},
\end{aligned} \right.
\end{equation}
where  we are using the notation $U_\alpha=-\Delta_{-\alpha}$.

We will use the fundamental solution $F_\alpha(x,y;t)$ to the Schr\"odinger equation (\ref{beta})  for obtaining the propagator (unitary group). Then we have the representation
\begin{equation}\label{uni}
e^{-itU_\alpha}f(x)=\int_{\mathbb R} F_\alpha(x,y;t)f(y)dy.
\end{equation}
Indeed,  from \cite{ABD} and \cite{KS} we have for  $S(x;t)$ denoting the free propagator in $\mathbb R$, i.e.
\begin{equation}\label{free}
S(x,t)= \frac{e^{-x^2/{4it}}}{(4i\pi t)^{1/2}},\quad t>0
\end{equation}
 and so $e^{it\Delta}f(x)=S(x;t)\ast_x f(x)$,  the following expression for $a\alpha \neq -1$:
\begin{enumerate}
\item[1)] For $\alpha>0$
\begin{equation}\label{ker51}
F_\alpha(x,y;t)=S(x-y;t)- \frac{1}{2\pi i}\int_{\mathbb R}e^{-i\xi ^2t}\frac{f_\alpha(x,y;\xi)}{(2\xi+i\alpha)^2+\alpha^2e^{i4\xi a}}d\xi
\end{equation}
with $f_\alpha(x,y;\xi)=\sum_{j=1}^4 L_\alpha^j(x,y;\xi)$ and

\begin{align}
L_\alpha^1(x,y;\xi)&=-\alpha(2\xi+i\alpha)e^{i\xi|x+a|}e^{i\xi|y+a|},\qquad L_\alpha^4(x,y;\xi)=L_\alpha^1(-x,-y;\xi)\\
L_\alpha^2(x,y;\xi)&=i\alpha^2e^{2i\xi a}e^{i\xi|x+a|}e^{i\xi|y-a|},\qquad\;\;\;\;\;\;\; L_\alpha^3(x,y;\xi)=L_\alpha^2(-x,-y;\xi).
\end{align}

\item[2)] For $\alpha<0$,
\begin{equation}\label{ker6}
F_\alpha(x,y;t)=e^{-it\gamma_1} \Gamma_1(x) \Gamma_1(y)+e^{-it\gamma_2} \Gamma_2(x) \Gamma_2(y) +F_{-\alpha}(x,y;t)
\end{equation}

where $ \Gamma_1$ and  $ \Gamma_2$ are the normalized eigenfunction associated with the eigenvalues $\gamma_1$ and $\gamma_2$.
\end{enumerate}

\begin{remark} We observe the following:
\begin{itemize}
\item[1) ]  The case $a\alpha=-1$ is assumed for technical reasons (see \cite{KS})
\item[2) ] For  $a\alpha=-1$ we obtain that only $\gamma_1$ remains as an eigenvalue in the discrete spectrum.
\end{itemize}
\end{remark}

\subsubsection{The case $\mu(x)=\beta \delta'$.}

Next we determine the linear propagator $J_\beta(t)=e^{i(\Delta+\beta\delta')t}$ (unitary group) determined by the linear system associated with (\ref{SCH0}),
\begin{equation}\label{beta}
\left \{
\begin{aligned}
i\partial_{t}u&=-(\Delta u+\beta\delta')u\equiv K_\beta u, \\
u(0)  &  =u_{0},
\end{aligned} \right.
\end{equation}
where  we are using the notation $K_\beta=-\Delta_{-\beta}$.

We will use the fundamental solution $S_\beta(x,y;t)$ to the Schr\"odinger equation (\ref{beta})  for obtaining the propagator (unitary group). Then we have the representation
\begin{equation}\label{uni}
e^{-itK_\beta}f(x)=\int_{\mathbb R} S_\beta(x,y;t)f(y)dy.
\end{equation}
Indeed,  from \cite{AGHH} we have for  $S(x;t)$ denoting the free propagator in (\ref{free}) the following:
\begin{enumerate}
\item[1)] For $\beta>0$
\begin{equation}\label{ker5}
S_\beta(x,y;t)=S(x-y;t)+\text{sgn}(xy) S(|x|+|y|;t)+ \frac{2}{\beta}\int_0^\infty \text{sgn}(xy)e^{-\frac{2}{\beta} s} S(s+|x|+|y|;t)ds
\end{equation}

\item[2)] For $\beta<0$,
\begin{equation}\label{ker6}
\begin{aligned}
S_\beta(x,y;t)=&S(x-y;t)+\text{sgn}(xy) S(|x|+|y|;t)+ e^{i\frac{4}{\beta^2} t} \Phi_\beta(x) \Phi_\beta(y)\\
&-\frac{2}{\beta} \int_0^\infty \text{sgn}(xy)e^{\frac{2}{\beta} s} S(s-|x|-|y|;t)ds
\end{aligned}
\end{equation}

where $ \Phi_\beta$ is defined in Theorem \ref{resol5d}.
\end{enumerate}

\subsubsection{Dispersive Estimates}

The following proposition extends the well known estimates for the free propagator $e^{it\Delta}$,
\begin{equation}\label{stric0}
\left\Vert e^{it\Delta}f(t)\right\Vert _{p^{\prime}}\leq C_0t^{-\frac{1}{2}%
(\frac{2}{p}-1)}\left\Vert f\right\Vert _{p},
\end{equation}
to the the groups $G_\sigma(t)=e^{-it H_\sigma}$, $M_\alpha(t)=e^{-it U_\alpha}$ and $J_\beta(t)=e^{i(\Delta+\beta\delta')t} $ in the  one-dimensional case. We denote by  $W_1$ the group $G_\sigma$,  $W_2$ the group $M_\alpha$ and by $W_3$ the group $J_\beta$ .

\begin{proposition} \label{strich}
Let $p\in [1, 2]$ and $p'$ be such that $\frac{1}{p}+\frac{1}{p'}=1$. Then we have:

Suppose that $u(x,t)=W_i(t)f(x)$, $i=1,2,3$, is the solution of the linear equation (\ref{SCH1}), (\ref{beta3}) and  (\ref{beta}), respectively. Then:
\begin{itemize}
\item[1)] for $\sigma\geqq 0$, $\alpha \geqq 0$ and $\beta\geqq 0$,
\begin{equation}\label{stric2}
\left\Vert W_i(t)f\right\Vert _{p^{\prime}}\leq C|t|^{-\frac{1}{2}%
(\frac{2}{p}-1)}\left\Vert f\right\Vert _{p},
\end{equation}
where $C>0$ is independent of $f$ and $t\neq0$.
\item[2)] for $\sigma< 0$, $a\leqq -\frac{1}{\alpha}$ ($\alpha<0$) and $\beta<0$,
\begin{equation}\label{stric3}
\left\Vert W_i(t)f-e^{i\alpha_i t} P_if\right\Vert _{p^{\prime}}\leq C|t|^{-\frac{1}{2}%
(\frac{2}{p}-1)}\left\Vert f\right\Vert _{p},
\end{equation}
where $\alpha_1=\frac{\sigma^2}{4}$, $P_1f= \langle f, \Psi_\sigma\rangle \Psi_\sigma$, $\alpha_2=-\gamma_1(a,\alpha)$, $P_2f= \langle f, \Gamma_1\rangle \Gamma_1$ and $\alpha_3=\frac{4}{\beta^2}$, $P_3f= \langle f, \Phi_\beta\rangle \Phi_\beta$, and $C>0$ is independent of $f$ and $t\neq0$.

\item[3)] for $\alpha<0$,  $a> -\frac{1}{\alpha}$,
\begin{equation}\label{stric4}
\left\Vert W_2(t)f-\sum_{j=1}^2e^{-i\gamma_j t} Q_jf\right\Vert _{p^{\prime}}\leq C|t|^{-\frac{1}{2}%
(\frac{2}{p}-1)}\left\Vert f\right\Vert _{p},
\end{equation}
where $\gamma_j=\gamma_j(a,\alpha)$, $j=1,2$,
$Q_1f=  \langle f, \Gamma_1\rangle \Gamma_1$, $Q_2f=  \langle f, \Gamma_2\rangle \Gamma_2$, and $C>0$ is independent of $f$ and $t\neq0$.

\end{itemize}
\end{proposition}

\textbf{Proof.} i) We consider $\sigma>0$. Initially  $G_{\sigma}(t)$ is a unitary group on $L^2(\Bbb R)$, $\|G_{\sigma}\phi(t)\|_2=\|\phi\|_2$ for all $t\in \mathbb R$. Let $\phi\in L^1(\mathbb R)$ and $R\phi(x)=\phi(-x)$. Then for $\phi^-=\phi \chi^0_-$ and $\phi^+=\phi R\chi^0_+$ we have the decomposition $\phi=\phi^-+R\phi^+$. Hence since $supp\; \phi^+\subset (-\infty, 0]$, $RG_{\sigma}=G_{\sigma}R$ and $R(f\ast Rg)=(Rf)\ast g$
we obtain from the following equality,
\begin{align}
G_{\sigma}\phi(t)&=[e^{it\Delta}\phi^-+e^{it\Delta}(\phi^-\ast \rho_\sigma)]\chi^0_-+e^{it\Delta}(\phi^-\ast \tau_\sigma)\chi^0_+\\
&+[e^{it\Delta}R\phi^+ +e^{it\Delta}(\phi^+\ast \rho_\sigma)]\chi^0_ + +e^{it\Delta}(R\phi^+\ast R\tau_\sigma)\chi^0_-.
\end{align}
 Therefore from (\ref{stric0}) and  applying Young's inequality we obtain for $t\neq 0$
\begin{align}
 \|G_{\sigma}\phi(t)\|_\infty&\leqq \frac{C_0}{\sqrt{|t|}}(\|\phi^-\|_1+\|\phi^-\ast \rho_\sigma\|_1+\|\phi^-\ast \tau_\sigma\|_1+\|R\phi^+\|_1+\|\phi^+\ast \rho_\sigma\|_1 +\|R\phi^+\ast R\tau_\sigma\|_1)\\
& \leqq \frac{C_0}{\sqrt{|t|}}(3+4\|\rho_\sigma\|_1)\|\phi\|_1=\frac{C}{\sqrt{|t|}}\|\phi\|_1
 \end{align}
 where $C=C(\sigma)$. By the Riesz-Thorin interpolation theorem  we obtain (\ref{stric2}). The case  $\sigma<0$ follows similarly from the expression (\ref{negpro}).

 ii) Let $\beta>0$. From (\ref{ker5}) we obtain immediate for $x,y\in \mathbb R$ that
 $$
 |S_\beta(x,y,t)|\leqq C|t|^{-\frac12}.
 $$
 So, from (\ref{uni}) we have $\|J_\beta(t)\phi\|_\infty\leqq \frac{C}{\sqrt{|t|}}\|\phi\|_1$. By following a similar analysis as in the later case we obtain (\ref{stric2}).

 iii) From the representations in (\ref{negpro}) and (\ref{ker6}) we get immediately (\ref{stric3}).

 iv) The case of the group $W_2$ for $a\alpha \neq -1$, it follows of the estimate
 $$
 \|e^{-it U_\alpha}P_c f\|_\infty\leqq C t^{-1/2}\|f\|_1,
 $$
 for $t>0$ and $P_c$ being the spectral projector of $U_\alpha$ on its  continuous spectrum.
 \fin

\section{\bigskip Weak-$L^{p}$ Solutions}

In this section we focus our study of global solutions for  the Cauchy problem
\begin{equation}
\left\{
\begin{aligned}
i\partial_{t}u+\Delta u+  \mu(x)u&  =\lambda\left\vert u\right\vert
^{\rho-1}u,\ x\in\mathbb{R},\ t\in\mathbb{R},\label{SCHd}\\
u(x,0)  &  =u_{0},
\end{aligned}\right.
\end{equation}
for $\mu(x)=\sigma \delta$ in the spaces $L^{(p, \infty)}(\mathbb R)$, which are called weak-$L^p$ or Marcinkiewicz spaces. The cases $\mu(x)=\alpha(\delta(\cdot-a)+ \delta(\cdot+a))$ and  $\mu(x)=\beta \delta'$ are treatment similarly.

We start by recalling some facts about the weak spaces $L^{(p, \infty)}(\mathbb R)$. For $1<r\leq
\infty,$ we recall that a measurable function $f$ defined on $\mathbb{R}$
belongs to $L^{(r,\infty)}(\mathbb{R})$ if the norm
\[
\Vert f\Vert_{(r,\infty)}=\sup_{t>0}t^{\frac{1}{r}}f^{\ast\ast}(t)
\]
is finite, where
\[
f^{\ast\ast}(t)=\frac{1}{t}\int_{0}^{t}f^{\ast}(s)\mbox{ }ds,
\]
and $ f^{\ast}$ is the {\it decreasing rearrangement} of $f$ with regard to the Lebesgue measure $\nu$, namely,
\[
f^{\ast}(t)=\inf\{s>0:\nu(\{x\in\mathbb{R}:|f(x)|>s\})\leq t\},\text{ }t>0,
\]
The space $L^{(r,\infty)}$ with the norm $\Vert f\Vert_{(r,\infty)}$ is a Banach space. We have the continuous inclusion $L^{r}(\mathbb{R})\subset$ $L^{(r,\infty)}(\mathbb{R})$. Moreover, the H\"{o}lder's inequality holds true in this framework, namely
\begin{equation}
\Vert fg\Vert_{(r,\infty)}\leq C\Vert f\Vert_{(q_{1},\infty)}\Vert
g\Vert_{(q_{2},\infty)}, \label{Holder1}%
\end{equation}
for $1<q_{1},q_{2}<\infty$, $\frac{1}{q_{1}}+\frac{1}{q_{2}}<1$ and $\frac
{1}{r}=\frac{1}{q_{1}}+\frac{1}{q_{2}}$, where $C>0$ depends only on $r$. Lastly, we have the Lorentz spaces $L^{(p,q)}(\mathbb R)$ that can be constructed via real interpolation; indeed, $L^{(p,q)}(\mathbb R)=(L^{1}(\mathbb R),L^{\infty}(\mathbb R))_{1-\frac{1}{p},q}$, $1<p<\infty$. They have the interpolation property
\begin{equation}\label{inter}
(L^{(p_0, q_0)}(\mathbb R), L^{(p_1, q_1)}(\mathbb R))_{\theta, q}=L^{(p,q)}(\mathbb R),
\end{equation}
provided $0<p_0<p_1<\infty$, $0<\theta<1$, $\frac{1}{p}=\frac{1-\theta}{p_0}+\frac{\theta}{p_1}$, $1\leqq q_0, q_1, q\leqq \infty$, where $(\cdot,\cdot)_{\theta, q}$ stands for the real interpolation spaces constructed via the $K$-method. For
further details about weak-$L^{r}$ and Lorentz spaces see \cite{BL} and Grafakos \cite{Gra}.

From (\ref{inter}) we obtain our main estimate for the group $G_\sigma$ in Lorentz spaces. A similar result is obtained for the groups $M_\alpha$ and $J_\beta$.

\begin{lemma}\label{grupint1}
Let $1\leqq d\leqq \infty$, $p' \in (2,\infty)$, and $p\in (1,2)$. If $p'$ satisfies $\frac{1}{p}+\frac{1}{p'}=1$, then there exists a constant $C>0$ such that:
\begin{enumerate}
\item[1)] for $\sigma\geqq 0$,
\begin{equation}\label{grupint2}
\left\Vert G_{\sigma}(t)f\right\Vert _{(p^{\prime},d)}\leq C|t|^{-\frac
{1}{2}(\frac{2}{p}-1)}\left\Vert f\right\Vert _{(p, d)},
\end{equation}
for all $f\in L^{(p,d)}(\mathbb R)$ and all $t\neq 0$.

\item[2)] for $\sigma< 0$,
\begin{equation}\label{grupint3}
\left\Vert G_{\sigma}(t)f-e^{i\frac{\sigma^2}{4}t} P_1f\right\Vert _{(p^{\prime},d)}\leq C|t|^{-\frac{1}{2}(\frac{2}{p}-1)}\left\Vert f\right\Vert _{(p, d)},
\end{equation}
for all $f\in L^{(p,d)}(\mathbb R)$ and all $t\neq 0$.
\end{enumerate}
\end{lemma}

\textbf{Proof.} We only consider the case $\sigma\geqq 0$ because for $\sigma<0$ the analysis is similar. Let fixed $t\neq 0$ and let $1<p_0<p<p_1<2$. From the $L^p=L^{(p,p)}$ estimate of the Schr\" odinger group in Proposition \ref{strich}, we have that $G_\sigma (t):L^{p_0}\to L^{p_0'}$ and $G_\sigma (t):L^{p_1}\to L^{p_1'}$ satisfy
$$
\|G_\sigma (t)\|_{p_0\to p_0'}\leqq C|t|^{-\frac
{1}{2}(\frac{2}{p_0}-1)},\qquad \|G_\sigma (t)\|_{p_1\to p_1'}\leqq C|t|^{-\frac
{1}{2}(\frac{2}{p_1}-1)}
$$
with $\frac{1}{p_0}+\frac{1}{p_0'}=1$ and $ \frac{1}{p_1}+\frac{1}{p_1'}=1$. Hence, for $\lambda\in (0,1)$, $\frac{1}{p}=\frac{1-\lambda}{p_0}+\frac{\lambda}{p_1}$, and $\frac{1}{p' }=\frac{1-\lambda}{p'_0}+\frac{\lambda}{p'_1}$, we obtain from (\ref{inter}) that
$$
\|G_\sigma (t)\|_{(p,d) \to (p',d)}\leqq \|G_\sigma (t)\|_{p_0\to p_0'}^\lambda\|G_\sigma (t)\|_{p_1\to p_1'}^{(1-\lambda)}\leqq C |t|^{-\frac{1}{2}(\frac{2}{p}-1)},
$$
which gives (\ref{grupint2}).

\fin

Next we establish the main results of this section.  From now on we focus in the case of $\mu(x)=\sigma\delta$ in (\ref{SCHd}), since for the case of the potential being the two symmetric deltas and the derivative of the Dirac- delta we have a similar analysis. We start by defining $\mathcal{L}_{\vartheta
}^{\infty}$ as the Banach space of all Bochner measurable functions
$u:\mathbb{R}\rightarrow$ $L^{(\rho+1,\infty)}$ endowed with the norm%

\begin{equation}
\Vert u\Vert_{\mathcal{L}_{\vartheta }^{\infty}}=\sup_{-\infty<t<\infty
}|t|^{\vartheta }\Vert u(t)\Vert_{(\rho+1,\infty)}, \label{Norma1}%
\end{equation}
where%
\begin{equation}
\vartheta =\frac{1}{\rho-1}-\frac{1}{2(\rho+1)}. \label{exponent1}%
\end{equation}
Let us also define the initial data space $\mathcal{E}_{0}$ as the set of all
$u\in\mathcal{S}^{\prime}(\mathbb{R})$ such that the norm
\[
\Vert u_{0}\Vert_{\mathcal{E}_{0}}=\sup_{-\infty<t<\infty}|t|^{\vartheta }\Vert
G_{\sigma}(t)u_{0}\Vert_{(\rho+1,\infty)}%
\]
is finite. Throughout this paper we stand for $\rho_{0}=\frac
{3+\sqrt{17}}{2}>1$ the positive root of the equation $\rho
^{2}-3\rho-2=0$.

 From Duhamel's principle, (\ref{SCHd}) is formally equivalent to
the integral equation
\begin{equation}
u(t)=G_{\sigma}(t)u_{0}-i\lambda\int_{0}^{t}G_{\sigma}(t-s)[|u(s)|^{\rho-1}u(s)]ds,
\label{int1}%
\end{equation}
where $G_{\sigma}(t)=e^{i(\Delta+\sigma\delta)t}$ is the group determined
by the linear system associated with { (\ref{SCHd})}.

\begin{definition} A mild solution of the initial value problem
(\ref{SCHd}) is a complex-valued function $u\in\mathcal{L}_{\vartheta }^{\infty}$ satisfying (\ref{int1}).
\end{definition}

Our main results of this section read as follows.

\begin{theorem}
\label{GlobalTheo} Let $\sigma\geqq 0$, $\rho_{0}<\rho<\infty$ and $u_{0}\in
\mathcal{E}_{0}.$ There is $\varepsilon>0$ such that {if }$\left\Vert
u_{0}\right\Vert _{\mathcal{E}_{0}}\leq\varepsilon${\ then the IVP
(\ref{SCHd}) has a unique global-in-time mild solution
$u\in\mathcal{L}_{\vartheta }^{\infty}$ satisfying }$\Vert u\Vert_{\mathcal{L}%
_{\vartheta }^{\infty}}\leq2\varepsilon.${ } Moreover, the data-solution map
$u_{0}\mapsto u$ from $\mathcal{E}_{0}$ into $\mathcal{L}_{\vartheta }^{\infty}$
is locally Lipschitz.
\end{theorem}

\begin{remark} The proof of Theorem \ref{GlobalTheo} is based in an argument of fixed point, so by using the implicit function theorem is not difficult to show that the data-solution map
$u_{0}\mapsto u$ from $\mathcal{E}_{0}$ into $\mathcal{L}_{\vartheta }^{\infty}$
is smooth.
\end{remark}

\begin{remark}
\label{rem-local}(Local-in-time solutions) Let $1<\rho<\rho_{0}$, $d_{0}=\frac{1}{2}(\frac{\rho
-1}{\rho+1}),$ and $d_{0}<\zeta<\frac{1}{\rho}$. For $0<T<\infty$, consider
the Banach space $\mathcal{L}_{\zeta}^{T}$ of all Bochner measurable functions
$u:(-T,T)\rightarrow$ $L^{(\rho+1,\infty)}$ endowed with the norm
\[
\Vert u\Vert_{\mathcal{L}_{\zeta}^{T}}=\sup_{-T<t<T}|t|^{\zeta}\Vert
u(\cdot,t)\Vert_{(\rho+1,\infty)}.
\]
A local-in-time existence result in $\mathcal{L}_{\zeta}^{T}$ could be proved
for (\ref{SCHd}) by considering $u_{0}\in L^{(\frac{\rho+1}{\rho},\infty)}%
(\mathbb{R})$ and small $T>0$ (see \cite{BFV}).
\end{remark}

In the sequel we give an asymptotic stability result for the obtained solutions.

\begin{theorem}
\label{TeoAssin}(Asymptotic Stability) Under the hypotheses of Theorem \ref{GlobalTheo} let $u$ and
$v$ be two solutions of (\ref{int1}) obtained through Theorem \ref{GlobalTheo}
with initial data ${u}_{0}$ and $v_{0},$ respectively. We have that%
\begin{equation}
\lim_{\left\vert t\right\vert \rightarrow\infty}\left\vert t\right\vert
^{\vartheta }\left\Vert u(\cdot,t)-v(\cdot,t)\right\Vert _{(\rho+1,\infty)}=0
\label{as1}%
\end{equation}
if only if
\begin{equation}
\lim_{\left\vert t\right\vert \rightarrow\infty}\left\vert t\right\vert
^{\vartheta }\left\Vert G_{\sigma}(t)(u_{0}-v_{0})\right\Vert _{(\rho+1,\infty)}=0
\label{as2}%
\end{equation}
The condition (\ref{as2}) holds, in particular, for $u_{0}-v_{0}\in
L^{(\frac{\rho+1}{\rho},\infty)}.$
\end{theorem}

\subsection{\bigskip Nonlinear Estimate}

In this subsection we give the nonlinear estimate essential in the proof of Theorem \ref{GlobalTheo}. We start by  recalling the Beta function
$$
B(\nu,\eta)=\int_{0}^{1}(1-s)^{\nu-1}%
s^{\eta-1}ds,
$$
 which is finite for all $\nu>0$ and $\eta>0.$ So, for $k_{1}%
,k_{2}<1$ and $t>0,$ the change of variable $s\rightarrow st$ yields
\begin{equation}
\int_{0}^{t}(t-s)^{-k_{1}}s^{-k_{2}}ds=t^{1-k_{1}-k_{2}}\int_{0}%
^{1}(1-s)^{-k_{1}}s^{-k_{2}}ds=t^{1-k_{1}-k_{2}}B(1-k_{1},1-k_{2})<\infty.
\label{Beta}%
\end{equation}

Next we denote by
\begin{equation}
\mathcal{N}(u)=-i\lambda \int_{0}^{t}G_{\sigma}(t-s)[\left\vert u(s)\right\vert ^{\rho-1
}u(s)]ds \label{non1}%
\end{equation}
the nonlinear part of the integral equation (\ref{int1}). We have the following estimate in order to apply a point fixed argument.

\begin{lemma}
\label{lem8}Let $\sigma\geqq 0$ and $\rho_{0}<\rho<\infty$. There is a
constant $K>0$ such that
\begin{equation}
\Vert\mathcal{N}(u)-\mathcal{N}(v)\Vert_{\mathcal{L}_{\vartheta }^{\infty}}\leq
K \Vert u-v\Vert_{\mathcal{L}_{\vartheta }^{\infty}}(\Vert u\Vert_{\mathcal{L}%
_{\vartheta }^{\infty}}^{\rho-1}+\Vert v\Vert_{\mathcal{L}_{\vartheta }^{\infty}%
}^{\rho-1}) \label{est-non1}%
\end{equation}
for all $u,v\in\mathcal{L}_{\vartheta }^{\infty}$.
\end{lemma}

\textbf{Proof.} Without loss of generality, we assume $t>0.$ It
follows from (\ref{grupint2}) with $p=\frac{{\rho+1}}{\rho}$, $d=\infty$ and
H\"older inequality (\ref{Holder1}) and $\||f|^r\|_{(p, \infty)}=\|f\|^r_{(rp, \infty)}$  that
\begin{align}
\Vert\mathcal{N}(u)-\mathcal{N}(v)\Vert_{(\rho+1,\infty)}  &  \leq\int_{0}%
^{t}\Vert G_{\sigma}(t-s)(\left\vert u\right\vert ^{\rho-1}u-\left\vert
v\right\vert ^{\rho-1}v)\Vert_{(\rho+1,\infty)}ds\nonumber\\
&  \hspace{-3.3cm}\leq C\int_{0}^{t}(t-s)^{-\frac{1}{2}(\frac{2\rho}{\rho
+1}-1)}\Vert(\left\vert u-v\right\vert )(\left\vert u\right\vert ^{\rho
-1}+\left\vert v\right\vert ^{\rho-1})\Vert_{(\frac{{\rho+1}}{\rho}{,\infty)}%
}ds\nonumber\\
&  \hspace{-3.3cm}\leq C\int_{0}^{t}(t-s)^{-\zeta}\Vert u-v\Vert
_{(\rho+1{,\infty)}}\left(  \Vert u\Vert_{(\rho+1{,\infty)}}^{\rho-1}+\Vert
v\Vert_{(\rho+1{,\infty)}}^{\rho-1}\right)  ds. \label{aux2}%
\end{align}
Next, notice that $\zeta=\frac{1(\rho-1)}{2(\rho+1)}<1$ and $\vartheta \rho<1$ when
$\rho_{0}<\rho<\infty.$ Thus, by using (\ref{Beta}), the r.h.s of
(\ref{aux2}) can be bounded by
\begin{align*}
&  \leq C\left(  \sup_{0<t<\infty}t^{\vartheta }\Vert u-v\Vert_{(\rho+1{,\infty)}%
}\sup_{0<t<\infty}\left(  t^{\vartheta (\rho-1)}\Vert u\Vert_{(\rho+1{,\infty)}%
}^{\rho-1}+t^{\vartheta (\rho-1)}\Vert v\Vert_{(\rho+1{,\infty)}}^{\rho-1}\right)
\right)  \times\int_{0}^{t}(t-s)^{-\zeta}s^{-\vartheta \rho}ds\\
&  =CB(1-\zeta, 1-\vartheta  \rho) t^{1-\vartheta \rho-\zeta}\left(  \Vert u-v\Vert_{\mathcal{L}_{\vartheta
}^{\infty}}(\Vert u\Vert_{\mathcal{L}_{\vartheta }^{\infty}}^{\rho-1}+\Vert
v\Vert_{\mathcal{L}_{\vartheta }^{\infty}}^{\rho-1})\right)  ,
\end{align*}
which implies (\ref{est-non1}), because $\zeta+\rho\vartheta =-\vartheta -1.$ \fin

\subsection{Proof of Theorem \ref{GlobalTheo} }

Consider the map $\Phi$ defined on $\mathcal{L}_{\vartheta }^{\infty}$ by
\begin{equation}
\Phi(u)=G_{\sigma}(t)u_{0}+\mathcal{N}(u)\label{map1}%
\end{equation}
where $\mathcal{N}(u)$ is given in (\ref{non1}). Let $\mathcal{B}
_{\varepsilon}=\{u\in\mathcal{L}_{\vartheta }^{\infty};\Vert u\Vert_{\mathcal{L}%
_{\vartheta }^{\infty}}\leq2\varepsilon\}$ where $\varepsilon>0$ will be chosen
later. Lemma \ref{lem8} implies that
\begin{align}\label{aux}
\Vert\Phi(u)-\Phi(v)\Vert_{\mathcal{L}_{\vartheta }^{\infty}} &  =\Vert
\mathcal{N}(u)-\mathcal{N}(v)\Vert_{\mathcal{L}_{\vartheta }^{\infty}}\nonumber\leq K  \Vert u-v\Vert_{\mathcal{L}_{\vartheta }^{\infty}}(\Vert u\Vert
_{\mathcal{L}_{\vartheta }^{\infty}}^{\rho-1}+\Vert v\Vert_{\mathcal{L}_{\vartheta
}^{\infty}}^{\rho-1})\\
&  \leq2^{\rho}\varepsilon^{\rho-1}K \Vert u-v\Vert_{\mathcal{L}_{\vartheta
}^{\infty}},
\end{align}
for all $u,v\in \mathcal{B}_{\varepsilon}.$ Since
\[
\Vert G_{\sigma}(t)u_{0}\Vert_{\mathcal{L}_{\vartheta }^{\infty}}=\Vert u_{0}%
\Vert_{\mathcal{E}_{0}}\leq\varepsilon,
\]
and by using  inequality (\ref{est-non1}) with $v=0$ we obtain%

\begin{align}
\Vert\Phi(u)\Vert_{\mathcal{L}_{\vartheta }^{\infty}}  &  \leq\Vert G_{\sigma
}(t)u_{0}\Vert_{\mathcal{L}_{\vartheta }^{\infty}}+\Vert\mathcal{N}(u)\Vert
_{\mathcal{L}_{\vartheta }^{\infty}}\nonumber  \leq\Vert G_{\sigma}(t)u_{0}\Vert_{\mathcal{L}_{\vartheta }^{\infty}}+K  \Vert
u\Vert_{\mathcal{L}_{\vartheta }^{\infty}}^{\rho}\nonumber\\
&  \leq\varepsilon+2^{\rho}\varepsilon^{\rho}K  \leq2\varepsilon, \label{aux6}%
\end{align}
provided that $2^{\rho}\varepsilon^{\rho-1}K <1$ and $u\in \mathcal{B}_{\varepsilon}.$ It follows that $\Phi:\mathcal{B}_{\varepsilon}\rightarrow \mathcal{B}_{\varepsilon}$ is a contraction, and then it has a fixed point $u\in \mathcal{B}_{\varepsilon},$ $\Phi(u)=u$, which is the unique solution for the integral equation (\ref{int1}) satisfying $\Vert u\Vert
_{\mathcal{\mathcal L}^\infty_{\vartheta }}\leq2\varepsilon.$

In view of (\ref{aux}), if $u,v$ are two integral solutions with respective
data $u_{0},v_{0}$, then
\begin{align*}
\Vert u-v\Vert_{\mathcal{L}_{\vartheta }^{\infty}}  &  =\Vert G_{\sigma}%
(t)(u_{0}-v_{0})\Vert_{\mathcal{L}_{\vartheta }^{\infty}}+\Vert\mathcal{N}%
(u)-\mathcal{N}(v)\Vert_{\mathcal{L}_{\vartheta }^{\infty}}\\
&  \leq\Vert u_{0}-v_{0}\Vert_{\mathcal{E}_{0}}+2^{\rho}\varepsilon^{\rho
-1}K\Vert u-v\Vert_{\mathcal{L}_{\vartheta }^{\infty}},
\end{align*}
which, due to $2^{\rho}\varepsilon^{\rho-1}K <1,$ yields the Lipschitz
continuity of the data-solution map \fin

\subsection{Proof of Theorem \ref{TeoAssin}}

We will only prove that (\ref{as2}) implies (\ref{as1}). The converse follows
similarly (in fact it is easier) and it is left to the reader. For that matter, we subtract the
integral equations verified by $u$ and $v$ in order to obtain
\begin{align}
t^{\vartheta }\left\Vert u(\cdot,t)-v(\cdot,t)\right\Vert _{(\rho+1,\infty)}  &
\leq t^{\vartheta }\left\Vert G_{\sigma}(t)(u_{0}-v_{0})\right\Vert _{(\rho
+1,\infty)}+\nonumber\\
&  +t^{\vartheta }\left\Vert \int_{0}^{t}G_{\sigma}(t-s)(\left\vert u\right\vert
^{\rho-1}u-\left\vert v\right\vert ^{\rho-1}v)ds\right\Vert _{(\rho+1,\infty
)}. \label{diff2}%
\end{align}
Since $\Vert u\Vert_{\mathcal{L}^\infty_{\vartheta }},\Vert v\Vert_{\mathcal{L}^\infty_{\vartheta }%
}\leq2\varepsilon,$ we can estimate the second term in R.H.S. of (\ref{diff2})
as follows.
\begin{align}
&  I(t)=t^{\vartheta }\left\Vert\int_{0}^{t}G_{\sigma}(t-s)[\left\vert u\right\vert
^{\rho-1}u-\left\vert v\right\vert ^{\rho-1}v]ds\right\Vert_{(\rho+1,\infty
)}\nonumber\\
&  \leq Ct^{\vartheta }\int_{0}^{t}(t-s)^{-\frac{1}{2}(\frac{2\rho}{\rho+1}%
-1)}s^{-\vartheta \rho}s^{\vartheta }\Vert u(\cdot,s)-v(\cdot,s)\Vert_{(\rho
+1,\infty)}ds\left(  \Vert u\Vert_{\mathcal{L}^\infty_{\vartheta }}^{\rho-1}+\Vert
v\Vert_{\mathcal{L}^\infty_{\vartheta }}^{\rho-1}\right) \nonumber\\
&  \leq C2^{\rho}\varepsilon^{\rho-1}t^{\vartheta }\int_{0}^{t}(t-s)^{-\zeta
}s^{-\vartheta \rho}s^{\vartheta }\Vert u(\cdot,s)-v(\cdot,s)\Vert_{(\rho+1,\infty
)}ds, \label{aux9}%
\end{align}
where $\zeta=\frac{1}{2}(\frac{2\rho}{\rho+1}-1).$ Now, recalling that
$\zeta+\vartheta \rho-\vartheta -1=-\vartheta ,$ the change of variable $s\longmapsto ts$
in (\ref{aux9}) leads us to
\begin{equation}
I(t)\leq C2^{\rho}\varepsilon^{\rho-1}\int_{0}^{1}(1-s)^{-\zeta}s^{-\vartheta
\rho}(ts)^{\vartheta }\Vert u(\cdot,ts)-v(\cdot,ts)\Vert_{(\rho+1,\infty)}ds.
\label{aux10}%
\end{equation}
Set
\begin{equation}
L=\limsup_{t\rightarrow\infty}t^{\vartheta }\Vert u(\cdot,t)-v(\cdot
,t)\Vert_{(\rho+1,\infty)}<\infty\label{defM}%
\end{equation}
and recall from proofs of Lemma \ref{lem8} and Theorem \ref{GlobalTheo} that
\[
K =C\int_{0}^{1}(1-s)^{-\zeta}s^{-\vartheta \rho
}ds\quad\text{ and }\quad 2^{\rho}\varepsilon^{\rho-1}K <1.
\]
Then, computing $\limsup_{t\rightarrow\infty}$ in (\ref{diff2}) and using
(\ref{aux10}), we get
\begin{align*}
L  &  \leq\left(  C2^{\rho}\varepsilon^{\rho-1}\int_{0}^{1}(1-s)^{-\zeta
}s^{-\vartheta \rho}ds\right)  L\\
&  =2^{\rho}\varepsilon^{\rho-1}KL
\end{align*}
and therefore $L=0$, as required. \fin

\section{\bigskip Existence of a invariant manifold of periodic orbits}

It is not clear for us whether the approach applied in the proof of Theorem \ref{GlobalTheo} for the case $\mu(x)=\sigma \delta$ in (\ref{SCHd}) with $\sigma\geqq 0$ can be applied for the case $\sigma<0$. Similar situation is happening for the cases $\mu(x)=\alpha(\delta(\cdot-a)+\delta(\cdot+a))$ and $\mu(x)=\beta \delta'$  with $\alpha<0$ and $\beta<0$, respectively.

But, for instance,  in the case $\sigma<0$ we can establish a nice qualitative behavior associated to the linear flow generated by equation  (\ref{SCHd}). In fact,   it follows from Theorem \ref{resol5b} that  the linear part of the NLS-$\delta$ equation  (\ref{SCHd}) has a two-dimensional manifold of periodic orbits, namely,
$$
E^p=\{\gamma e^{i\theta} \Phi_\sigma(x):\gamma\geqq 0\;\;\text{and}\;\;\theta\in [0, 2\pi]\}.
$$
So, the estimate (\ref{grupint3}) will imply immediately that  all solutions $u(t)$ of (\ref{SCHd}) with $\lambda=0$  and with initial conditions $u_0\in L^{(p,d)}(\mathbb R)$ will approach to one of the periodic orbits $\gamma e^{i(\frac{\sigma^2}{4}t+\theta)}\Phi_\sigma\in E^p$. More exactly,  we have   the following theorem.

\begin{theorem} \label{manifold} Let $\sigma< 0$. For $ d\in [1,\infty]$, $p'\in [1, \infty]$ and $p\in[1,2]$, we have that for $p'$ satisfying $\frac{1}{p}+\frac{1}{p'}=1$, the solution $u(t)$ of the linear equation associated to (\ref{SCHd}) with initial data $u(0)=u_0\in L^{(p,d)}(\mathbb R)$ satisfies
$$
\lim_{t\to \pm \infty} \|u(t)-\gamma_0 e^{i(\frac{\sigma^2}{4}t+\theta)}\Phi_\sigma\|_{(p',d)}=0,
$$
for $\gamma_0=|\langle u_0, \Phi_\sigma\rangle|$ and some $\theta\in [0,2\pi]$.
\end{theorem}

\begin{remark} 1) A similar result to that in Theorem \ref{manifold} can be obtained   for the linear equation associated to (\ref{SCHd}) in  the case of $\mu(x)=\beta \delta'$ with $\beta<0$ and for the linear equation  (\ref{beta3}) in  the case of $\mu(x)=\alpha(\delta(\cdot-a)+\delta(\cdot+a))$ with $a\leqq -\frac{1}{\alpha}$ and $\alpha<0$.

2) Note that $\gamma_0<\infty$. In fact, it is not difficult to see that for $
 \Psi_\sigma(x)=\sqrt{\frac{-\sigma}{2}}e^{\frac{\sigma}{2}|x|}$ we have for $s\geqq 0$ that $\Psi^*_\sigma(s)=\sqrt{\frac{-\sigma}{2}}e^{\frac{\sigma}{4}s}$. So for all $p, q\in (0,\infty)$ we obtain that
 $$
 \| \Psi_\sigma\|^q_{L^{(p,q)}}=\int_0^\infty \Big(t^{\frac{1}{p}}\Psi^*_\sigma(t)\Big)^q\frac{dt}{t}= \Big(\frac{-\sigma}{2}\Big)^{\frac{q}{2}}\Big(\frac{-4}{q\sigma}\Big)^{\frac{q}{p}} \Gamma \Big(\frac{q}{p}\Big),
 $$
 where $\Gamma$ represents the Gamma function. The case $q=\infty$ is immediate. Next, by the Hardy-Littlewood inequality for decreasing rearrangements and the H\"older inequality in the classical $L^p(d\nu)$ spaces, we obtain for $p\in [1,2]$, $p'$ such that $\frac{1}{p}+\frac{1}{p'}=1$ and for $r$ such that $\frac{1}{d}+\frac{1}{r}=1$, with $d\geqq 1$, the estimate

 \begin{align*}
 \beta_0\leqq \int_{\mathbb R} |u_0(x)|| \Psi_\sigma (x)|dx&\leqq \int_0^\infty u^*_0(t) \Psi^*_\sigma (t)dt=\int_0^\infty t^{\frac{1}{p}}u^*_0(t) t^{\frac{1}{p'}}\Psi^*_\sigma (t)\frac {dt}{t}\\
 &\leqq \|u_0\|_{(p,d)}\|\Psi_\sigma\|_{(p',r)}<\infty.
 \end{align*}

\end{remark}

\section{\bigskip Spaces based on Fourier transform}

\ In this section we consider the nonlinear Schr\"odinger equation
\begin{equation}\label{SCH-F1}
\left\{
\begin{aligned}
i\partial_{t}u+\Delta u+\mu(x)u  &  =\lambda u^\rho,\ x\in\Omega,\ t\in \mathbb{R},\\
u(x,0)  &  =u_{0},\text{ }x\in\Omega%
\end{aligned}\right.
\end{equation}
where $\mu\in
BC(\mathbb{R}^{n})$ (the space of all bounded continuous functions on $\mathbb{R}^{n})$, $\lambda=\pm1$ \ and $\rho\in\mathbb{N}$. Here we will consider
$\Omega=\mathbb{T}^{n}$ and $\Omega=\mathbb{R}^{n}$, i.e. the periodic and
nonperiodic cases, respectively. The nonlinearity $\lambda\left\vert u\right\vert ^{\rho-1}u$ could be considered in (\ref{SCH-F1}), however we prefer $u^{\rho}$ for the sake of simplicity of
the exposition. For more details, see Remark \ref{rem-Fourier2} below.

We start by defining the spaces for the nonperiodic case. We recall that if $\mathcal{M}(\mathbb{R}^{n})$ denotes the space of complex Radon measures on $\mathbb{R}^{n}$, then it is a  vector space and for $\nu\in\mathcal{M}(\mathbb{R}^{n})$, $\|\nu\|_{\mathcal{M}}=|\nu|(\mathbb{R}^{n})$ is a norm on it, where $|\nu|$ is the total variation of $\nu$ (we note that every measure in $\mathcal{M}(\mathbb{R}^{n})$ is automatically a finite Radon measure. Moreover, we can embed $L^1(\mathbb{R}^{n}, dm)$ into $\mathcal{M}(\mathbb{R}^{n})$ by identifying $f\in L^1(\mathbb{R}^{n}, dm)$ with the complex measure $d\nu= fdm$, and $\|\nu\|_{\mathcal{M}}=\int |f|dm$.  Next, every $\nu\in \mathcal{M}(\mathbb{R}^{n})$ defines a tempered distribution by $T_\nu(\varphi)=\int_{\mathbb{R}^{n}} \varphi(x)d\nu$, thereby identifying $ \mathcal{M}(\mathbb{R}^{n})$ with a subspace of $\mathcal S'$.

The Fourier transform on $L^1(\mathbb{R}^{n})$ can be extended of a natural form to $\mathcal{M}(\mathbb{R}^{n})$; if $\nu\in \mathcal{M}(\mathbb{R}^{n})$, the Fourier transform of $\nu$ is the function $\widehat{\nu}$ defined by
\begin{equation}\label{transf}
\widehat{\nu}(\xi)=\int e^{-2\pi i \xi\cdot x} d\nu(x),\quad \xi\in \mathbb{R}^{n}.
\end{equation}
Using that $ e^{-2\pi i \xi\cdot x} $ is uniformly continuous in $x$, it is not difficult to check that $\widehat{\nu} \in BC(\mathbb{R}^{n})$ and that $\|\widehat{\nu}\|_{\infty}\leqq \|\nu\|_{\mathcal{M}}$ (see Folland \cite{FO}). Similarly, we define the inverse Fourier transform $\check{\nu}$ of $\nu$ by $\check{\nu}(\xi)=\widehat{\nu}(-\xi)$ for $\xi\in \mathbb{R}^{n}$. Moreover, if $\mathcal F$ represents the Fourier transform  on $\mathcal S'$, then for every $\nu\in \mathcal{M}(\mathbb{R}^{n})$ we have $\mathcal F (\nu)=\widehat{\nu}$. Similarly, $\mathcal F^{-1} (\nu)=\check{\nu}$.

Recall that the space $\mathcal{M}(\mathbb{R}^{n})$ can be identified with a subspace of $\mathcal S'$. Hence, if we assume that $ \mu\in \mathcal S'$ and $\mathcal F(\mu)$ is a finite Radon measure, then for $\nu=\mathcal F(\mu)$ we  have
$$
\mu=\mathcal F^{-1} \mathcal F(\mu)=\mathcal F^{-1} (\nu)=\check{\nu}\in BC(\mathbb{R}^{n}).
$$

Next, by using the above identification between $\mathcal F$ and\;  $\widehat{}$\; \;,  we define the Banach space
\begin{equation}
\mathcal{I=}[\mathcal{M}(\mathbb{R}^{n})]^{\vee}=\{f\in\mathcal{S}^{\prime
}(\mathbb{R}^{n}):\widehat{f}\in\mathcal{M}(\mathbb{R}^{n})\}\subset
BC(\mathbb{R}^{n}), \label{esp2}%
\end{equation}
with norm
\begin{equation}
\left\Vert f\right\Vert _{\mathcal{I}}=\Vert\widehat{f}\;\Vert
_{\mathcal{M}}.\text{ } \label{norm1}%
\end{equation}

We note that in general $\mu\in \mathcal{I}$ may not to belong to $L^{p}(\mathbb{R}^{n})$, nor to $L^{p,\infty}(\mathbb{R}^{n}),$ with $p\neq\infty.$ In particular, $\mu\in \mathcal{I}$ may have infinite $L^{2}$-mass; for instance, if $\mu\equiv1$ then $\mathcal F(\mu)=\delta_{0}\in \mathcal{M}(\mathbb{R}^{n})$.

In the following we will consider $\mu, u_0\in \mathcal{I}$. The Cauchy problem (\ref{SCH-F1}) is formally equivalent to the
functional equation
\begin{equation}
u(t)=S(t)u_{0}+B(u)+L_{\mu}(u), \label{int2}
\end{equation}
where $S(t)=e^{it\Delta}$ is the Schr\"odinger group in $\mathbb{R}^{n}$, and the operators $L_{\mu}(u), B(u)$ are defined via Fourier transform by
\begin{equation}
\widehat{L_{\mu}(u)}=\int_{0}^{t}e^{-i\left\vert \xi\right\vert ^{2}%
(t-s)}(\widehat{\mu}\ast\widehat{u})(\xi,s)ds,\label{op1}%
\end{equation}
and
\begin{equation}
\widehat{B(u)}=\lambda\int_{0}^{t}e^{-i\left\vert \xi\right\vert ^{2}(t-s)}(\underbrace{\widehat{u}\ast\widehat{u}\ast...\ast\widehat{u}%
}_{\rho-times})
(\xi,s)ds, \label{op3}%
\end{equation}
for $\mu\in \mathcal{I}$ and $u\in L^{\infty}((-T,T);\mathcal{I})$. We recall that for arbitrary $\mu, \nu\in \mathcal{M}(\mathbb{R}^{n})$ their convolution $\widehat{\mu}\ast\widehat{\nu}\in \mathcal{M}(\mathbb{R}^{n})$ is defined by
$$
\widehat{\mu}\ast\widehat{\nu}(E)=\int\int\chi_E(x+y)d\mu(x)d\nu(y),
$$
for every Borel set $E$.

Let $\mathbb{T}^{n}=\mathbb{R}^{n}/{\mathbb{Z}^{n}}$ stand for the $n$-torus. We say that functions on $\mathbb{T}^{n}$ are functions $f:\mathbb{R}^{n}\to \mathbb{C}$ that satisfy $f(x+m)=f(x)$ for all $x\in \mathbb{R}^{n}$ and $m\in \mathbb{Z}^{n}$, which are called $1$-periodic in every coordinate. Let $\mathcal P=C^\infty_{per}=\{f:\mathbb{R}^{n}\to \mathbb{C}: f \;\text{is}\; C^\infty \;\text{and periodic with period}\; 1\}$. So $\mathcal{D}^{\prime}(\mathbb{T}^{n})$ is the set of all periodic distributions on $\mathcal P$. We say that $\mathcal T:\mathcal P\to \mathbb{C}$ is a {\it periodic distribution} if there exists a sequence $(\Psi_j)_{j\geqq 1}\subset \mathcal P$ such that
$$
\mathcal T(f)=\langle \mathcal T, f\rangle=\lim_{j\to \infty}\int_{[-1/2,1/2]^n} \Psi_j(x)f(x)dx,\qquad \; f\in \mathcal P.
$$
Above we have identify $\mathbb{T}^{n}$ with $[-1/2,1/2]^n$. For a complex-valued function $f\in L^1(\mathbb{T}^{n})$ and $m\in  \mathbb{Z}^{n}$, we define
$$
\widehat{f}(m)=\int_{[-1/2,1/2]^n}f(x)e^{-2\pi i x\cdot m}dx.
$$
We call $\widehat{f}(m)$ the $m$-th Fourier coefficient of $f$. The Fourier series of $f$ at $x\in \mathbb{T}^{n}$ is the series
$$
\sum_{m\in \mathbb{Z}^{n}} \widehat{f}(m) e^{2\pi i x\cdot m}.
$$
The Fourier transform of $\mathcal T\in \mathcal{D}^{\prime}(\mathbb{T}^{n})$ is the function $\widehat{\mathcal T}: \mathbb{Z}^{n}\to \mathbb C$ defined by the formula
$$
\widehat{\mathcal T}(m)=\langle \mathcal T, e^{-2\pi i  x\cdot m}\rangle,\quad m\in \mathbb{Z}^{n}.
$$

In the periodic case, we are going to study (\ref{SCH-F1}) in the space $\mathcal{I}_{per}$ which is defined by
\begin{equation}
\mathcal{I}_{per}=\{f\in\mathcal{D}^{\prime}(\mathbb{T}^{n}):\widehat{f}\in
l^{1}(\mathbb{Z}^{n})\} \label{peri1}%
\end{equation}
endowed with the norm
\begin{equation}
\left\Vert f\right\Vert _{\mathcal{I}_{per}}=\Vert\widehat{f}\;\Vert
_{l^{1}(\mathbb{Z}^{n})}\text{.} \label{peri2}%
\end{equation}
Here the IVP (\ref{SCH-F1}) is formally converted to
\begin{equation}
u(t)=S_{per}(t)u_{0}+B_{per}(u)+L_{\mu,per}(u), \label{int-per}%
\end{equation}
where, similarly to above, we define the operators in (\ref{int-per}) via
Fourier transform in $\mathcal{D}^{\prime}(\mathbb{T}^{n})$. Precisely, $S_{per}(t)$ is the Schrodinger group in $\mathbb{T}^{n}$
\begin{equation}
S_{per}(t)u_{0}=\sum_{m\in\mathbb{Z}^{n}}\widehat{u}_{0}(m)e^{-4\pi
^{2}i\left\vert m\right\vert ^{2}t}e^{2\pi ix\cdot m}, \label{sch-per}%
\end{equation}%

\begin{equation}
\widehat{L_{\mu,per}(u)}(m,t)=\int_{0}^{t}e^{-4\pi^{2}i\left\vert m\right\vert
^{2}(t-s)}(\widehat{\mu}\ast\widehat{u})(m,s)ds \label{Op-per1}%
\end{equation}
and
\begin{equation}
\widehat{B_{per}(u)}(m,t)=\int_{0}^{t}e^{-4\pi^{2}i\left\vert m\right\vert
^{2}(t-s)}(\underbrace{\widehat{u}\ast\widehat{u}\ast...\ast\widehat{u}%
}_{\rho-times})(m,s)ds,
\label{Op-per2}%
\end{equation}
for $u_{0},\mu\in \mathcal{I}_{per}$ and $u\in L^{\infty}((-T,T);\mathcal{I}_{per})$, where now the symbol $\ast$ denotes the discrete convolution
\[
\widehat{f}\ast\widehat{g}(m)=\sum_{\xi\in\mathbb{Z}^{n}}\widehat{f}%
(m-\xi)\widehat{g}(\xi).
\]

Throughout this section, solutions of (\ref{int2}) or (\ref{int-per}) will be
called mild solutions for the IVP (\ref{SCH-F1}), according the
respective case.

In the above framework, our local-in-time well-posedness result reads as follows.

\begin{theorem}
\label{teoF1}Let $1\leq\rho<\infty.$
\end{theorem}

\begin{itemize}
\item[(1)] (Periodic case) Let $u_{0}\in\mathcal{I}_{per}$ and $\mu$
$\in\mathcal{I}_{per}$. There is $T>0$ such that the IVP (\ref{SCH-F1}) has a unique mild solution $u\in L^{\infty}((-T,T);\mathcal{I}_{per})$ satisfying
\[
\sup_{t\in(-T,T)}\left\Vert u(\cdot,t)\right\Vert _{\mathcal{I}_{per}}%
\leq 2\left\Vert u_{0}\right\Vert _{\mathcal{I}_{per}}.
\]
Moreover, the data-map solution $u_{0}\rightarrow u$ is locally Lipschitz continuous
from $\mathcal{I}_{per}$ to {$L^{\infty}((-T,T);$}$\mathcal{I}_{per}\mathcal{)}.$

\item[(2)] (Nonperiodic case) Let $u_{0}\in\mathcal{I}$ and $\mu$
$\in\mathcal{I}$. The same conclusion of item (1) holds true by
replacing $\mathcal{I}_{per}$ by $\mathcal{I}$.
\end{itemize}

\begin{remark}
\label{rem-Fourier}In item (2) of the above theorem, one can show that the
solution $u(x,t)$ verifies $\widehat{u}(\xi,t)\in L^{1}(\mathbb{R}^{n}),$ for
all $t\in(-T,T),$ provided that $\widehat{u}_{0}\in L^{1}(\mathbb{R}^{n})$ and
$\widehat{\mu}\in L^{1}(\mathbb{R}^{n}).$ Then, due to Riemann-Lebesgue lemma,
it follows that
\[
u(x,t)\rightarrow0\text{ as }\left\vert x\right\vert \rightarrow\infty,\text{
for each }t\in(-T,T).
\]
\end{remark}

\subsection{\bigskip Nonlinear Estimate}

Before proceeding with the proof of Theorem \ref{teoF1}, let us recall the Young inequality for measures and
discrete convolutions (see Folland \cite{FO} and Iorio\&Iorio \cite{Ior}). For $\mu, \nu \in \mathcal{M}(\mathbb R^n)$ and $f, g\in l^{1}=l^{1}(\mathbb Z^n)$, we have the respective estimates
\begin{align}
\left\Vert \mu\ast \nu\right\Vert _{\mathcal{M}}  &  \leq\left\Vert \mu\right\Vert
_{\mathcal{M}}\left\Vert \nu\right\Vert _{\mathcal{M}}\label{Young}\\
\left\Vert f\ast g\right\Vert _{l^{1}}  &  \leq\left\Vert f\right\Vert
_{l^{1}}\left\Vert g\right\Vert _{l^{1}}. \label{Young-per}%
\end{align}

\begin{lemma}
\label{lem-est1}Let $1\leq\rho<\infty$ and $0<T<\infty.$

\begin{itemize}
\item[(i)] There exists a positive constant $K>0$ such that
\begin{align}
\sup_{t\in(-T,T)}\left\Vert S_{per}(t)u_{0}\right\Vert _{\mathcal{I}_{per}}
&  \leq\left\Vert u_{0}\right\Vert _{\mathcal{I}_{per}}\label{est1}\\
\sup_{t\in(-T,T)}\left\Vert L_{\mu,per}(u-v)\right\Vert _{\mathcal{I}_{per}}
&  \leq T\left\Vert \mu\right\Vert _{\mathcal{I}_{per}}\sup_{t\in
(-T,T)}\left\Vert u(\cdot,t)-v(\cdot,t)\right\Vert _{\mathcal{I}_{per}%
}\label{est2}\\
\sup_{t\in(-T,T)}\left\Vert B_{per}(u)-B_{per}(v)\right\Vert _{\mathcal{I}%
_{per}}  &  \leq KT\sup_{t\in(-T,T)}\left\Vert u(\cdot,t)-v(\cdot
,t)\right\Vert _{\mathcal{I}_{per}}\label{est3}\\
&  \times\left(  \sup_{t\in(-T,T)}\left\Vert u(\cdot,t)\right\Vert
_{\mathcal{I}_{per}}^{\rho-1}+\sup_{t\in(-T,T)}\left\Vert v(\cdot
,t)\right\Vert _{\mathcal{I}_{per}}^{\rho-1}\right)  ,\nonumber
\end{align}
for all $u_{0},\mu\in\mathcal{I}_{per}$ and $u,v\in${$L^{\infty}((-T,T);$}%
$\mathcal{I}_{per}).$

\item[(ii)] The above estimates still hold true with $S(t),L_{\mu}(u),B(u)$
and $\mathcal{I}$ in place of $S_{per}(t),L_{\mu,per}(t),B_{per}(u)$ and
$\mathcal{I}_{per}$, respectively.
\end{itemize}
\end{lemma}

\bigskip

\textbf{Proof.} \ We will only prove the item (i) because (ii) follows similarly
by using (\ref{Young}) instead of (\ref{Young-per}). From definition of
$S_{per}(t),$ we have that%
\[
\sup_{t\in(-T,T)}\left\Vert S_{per}(t)u_{0}\right\Vert _{\mathcal{I}_{per}%
}=\left\Vert \left(  \widehat{u}_{0}(m)e^{-4\pi^{2}i\left\vert m\right\vert
^{2}t}\right)  _{m\in\mathbb{Z}^{n}}\right\Vert _{l^{1}(\mathbb{Z}^{n})}%
\leq\left\Vert \hat{u}_{0}\right\Vert _{l^{1}(\mathbb{Z}^{n})}.
\]

The operator $L_{\mu,per}$ can be estimated as
\begin{align*}
\left\Vert L_{\mu,per}(u)\right\Vert _{\mathcal{I}_{per}}  &  =\left\Vert
\widehat{L_{\mu,per}(u)}\right\Vert _{l^{1}(\mathbb{Z}^{n})}\\
&  \leq\left\Vert \left(  \int_{0}^{t}e^{-4\pi^{2}i\left\vert m\right\vert
^{2}(t-s)}(\widehat{\mu}\ast\widehat{u})(m,s)ds\right)  _{m\in\mathbb{Z}^{n}%
}\right\Vert _{l_{1}(\mathbb{Z}^{n})}\\
&  \leq\sum_{m\in\mathbb{Z}^{n}}\left\vert \int_{0}^{t}e^{-4\pi^{2}i\left\vert
m\right\vert ^{2}(t-s)}(\widehat{\mu}\ast\widehat{u})(m,s)ds\right\vert \\
&  \leq\int_{0}^{t}\sum_{m\in\mathbb{Z}^{n}}\left\vert (\widehat{\mu}%
\ast\widehat{u})(m,s)\right\vert ds\\
&  \leq\int_{0}^{t}\left\Vert \widehat{\mu}\right\Vert _{l^{1}(\mathbb{Z}%
^{n})}\left\Vert \widehat{u}(\cdot,s)\right\Vert _{l^{1}(\mathbb{Z}^{n})}ds=\left\Vert \mu\right\Vert _{\mathcal{I}_{per}}\left\Vert u\right\Vert
_{L^{1}(0,T;\mathcal{I}_{per})} \leq T\left\Vert \mu\right\Vert _{\mathcal{I}_{per}}\left\Vert u\right\Vert
_{L^{\infty}(0,T;\mathcal{I}_{per})}.
\end{align*}
By elementary convolution properties and Young inequality (\ref{Young-per}),
it follows that
\begin{align*}
&  \left\Vert (\underbrace{\widehat{u}\ast\widehat{u}\ast...\ast\widehat{u}%
)}_{\rho-times}-\underbrace{(\widehat{v}\ast\widehat{v}\ast...\ast\widehat{v}%
)}_{\rho-times}\right\Vert _{l^{1}(\mathbb{Z}^{n})}\\
&  \leq\left\Vert \lbrack(\hat{u}-\hat{v})\ast\widehat{u}\ast...\ast
\widehat{u}+\hat{v}\ast(\widehat{u}-\widehat{v})\ast...\ast\widehat{u}+\hat
{v}\ast\widehat{v}\ast(\widehat{u}-\widehat{v})\ast...\ast\widehat{u}%
+...+\hat{v}\ast\widehat{v}\ast...\ast(\widehat{u}-\widehat{v})\right\Vert
_{l^{1}(\mathbb{Z}^{n})}\\
&  \leq\left\Vert (\hat{u}-\hat{v})\right\Vert _{l^{1}}\left\Vert \widehat
{u}\right\Vert _{l^{1}}^{\rho-1}+\left\Vert (\hat{u}-\hat{v})\right\Vert
_{l^{1}}\left\Vert \widehat{u}\right\Vert _{l^{1}}^{\rho-2}\left\Vert
\widehat{v}\right\Vert _{l^{1}}+...+\left\Vert (\hat{u}-\hat{v})\right\Vert
_{l^{1}}\left\Vert \widehat{u}\right\Vert _{l^{1}}\left\Vert \widehat
{v}\right\Vert _{l^{1}}^{\rho-2}+\left\Vert (\hat{u}-\hat{v})\right\Vert
_{l^{1}}\left\Vert \widehat{v}\right\Vert _{l^{1}}^{\rho-1}\\
&  \leq K\left\Vert (\hat{u}-\hat{v})\right\Vert _{l^{1}}\left(  \left\Vert
\widehat{u}\right\Vert _{l^{1}}^{\rho-1}+\left\Vert \widehat{v}\right\Vert
_{l^{1}}^{\rho-1}\right)
\end{align*}
Therefore%
\begin{align*}
\left\Vert B_{per}(u)(t)-B_{per}(v)(t)\right\Vert _{\mathcal{I}_{per}}  &  =\left\Vert
\widehat{B_{per}(u)}-\widehat{B_{per}(v)}\right\Vert _{l^{1}}\\
&  \leq\left\Vert \int_{0}^{t}e^{-4\pi^{2}i\left\vert \xi\right\vert
^{2}(t-s)}\left[  (\underbrace{\widehat{u}\ast\widehat{u}\ast...\ast
\widehat{u})}_{\rho-times}-\underbrace{(\widehat{v}\ast\widehat{v}\ast
...\ast\widehat{v})}_{\rho-times}\right]  ds\right\Vert _{l^{1}}\\
&  \leq\int_{0}^{t}\left\Vert \left[  (\underbrace{\widehat{u}\ast\widehat
{u}\ast...\ast\widehat{u})}_{\rho-times}-\underbrace{(\widehat{v}\ast\widehat
{v}\ast...\ast\widehat{v})}_{\rho-times}\right]  \right\Vert _{l^{1}}ds\\
&  \leq K\int_{0}^{t}\left\Vert \hat{u}-\hat{v}\right\Vert _{l^{1}}\left(
\left\Vert \widehat{u}\right\Vert _{l^{1}}^{\rho-1}+\left\Vert \widehat
{v}\right\Vert _{l^{1}}^{\rho-1}\right)  ds\\
&  \leq KT\left\Vert u-v\right\Vert _{L^{\infty}(0,T;\mathcal{I}_{per
}\mathcal{)}}\left(  \left\Vert u\right\Vert _{L^{\infty}(0,T;\mathcal{I}%
_{per}\mathcal{)}}^{\rho-1}+\left\Vert v\right\Vert _{L^{\infty}%
(0,T;\mathcal{I}_{per}\mathcal{)}}^{\rho-1}\right)  ,
\end{align*}
as required. \fin

\begin{remark}
\label{rem-Fourier2}
The approach employed here could be used to treat (\ref{SCH-F1}) with
the nonlinearity $\left\vert u\right\vert ^{\rho-1}u$ instead of $u^{\rho}.$
For $\rho$ odd, it would be enough to write   $\left\vert u\right\vert
^{\rho-1}u$ (in the above proof) as
\begin{align*}
\lbrack\left(  \left\vert u\right\vert ^{2}\right)  ^{\frac{\rho-1}{2}%
}u]^{\wedge}  & =[(u\cdot\overline{u})^{\frac{\rho-1}{2}}u]^{\wedge}\\
& =(\underbrace{\widehat{u}\ast\widehat{u}\ast...\ast\widehat{u})}_{\frac
{\rho-1}{2}-times}\ast(\underbrace{\widehat{\overline{u}}\ast\widehat
{\overline{u}}\ast...\ast\widehat{\overline{u}})}_{\frac{\rho-1}{2}-times}\ast
u
\end{align*}
and to note that $\widehat{\overline{u}}(\xi)=\overline{\widehat{u}}(-\xi)$
and $\left\Vert \overline{u}(\xi)\right\Vert _{\mathcal{I}_{per}}=\left\Vert u(\xi)\right\Vert _{\mathcal{I}_{per}}$.
\end{remark}

\subsection{Proof of Theorem \ref{teoF1}}

\smallskip\textbf{Proof of (1).} \noindent Consider the ball $\mathcal{B}%
_{\varepsilon}=\{u\in L^{\infty}(-T,T;\mathcal{I}_{per});\Vert u\Vert_{L^{\infty
}(-T,T;\mathcal{I}_{per})}\leq2\varepsilon\}$ endowed with the complete metric
$Z(\cdot,\cdot)$ defined by
$$
Z(u,v)=\Vert u-v\Vert_{L^{\infty}(-T,T;\mathcal{I}%
_{per})}.
$$
Let $\varepsilon=\Vert u_{0}\Vert_{\mathcal{I}_{per}}$ and $T>0$
such that
\begin{equation}
T(2\left\Vert \mu\right\Vert _{\mathcal{I}_{per}}+2^{\rho}\varepsilon^{\rho
-1}K)<1. \label{cond-small}
\end{equation}
Notice that $\varepsilon$ can be large. We shall show that the map
\begin{equation}
\Phi(u)=S_{per}(t)u_{0}+L_{\mu,per}(u)+B_{per}(u)\nonumber
\end{equation}
\newline is a contraction on $(\mathcal{B}_{\varepsilon},Z).$ Lemma
\ref{lem-est1} with $v=0$ yields%

\begin{align}
\Vert\Phi(u)\Vert_{L^{\infty}(-T,T;\mathcal{I}_{per})}  &  \leq\Vert
S_{per}(t)u_{0}\Vert_{L^{\infty}(-T,T;\mathcal{I}_{per})}+\Vert L_{\mu,per}(u)\Vert_{L^{\infty}(-T,T;\mathcal{I}_{per})}\nonumber\\
&+\left\Vert B_{per}(u)\right\Vert
_{L^{\infty}(-T,T;\mathcal{I}_{per})}\nonumber\\
&  \leq\Vert u_{0}\Vert_{\mathcal{I}_{per}}+T\left\Vert \mu\right\Vert
_{\mathcal{I}_{per}}\Vert u\Vert_{L^{\infty}(-T,T;\mathcal{I}_{per})}+TK\Vert
u\Vert_{L^{\infty}(-T,T;\mathcal{I}_{per})}^{\rho}\nonumber\\
&  \leq\varepsilon+2\varepsilon T\left\Vert \mu\right\Vert _{\mathcal{I}%
_{per}}+2^{\rho}\varepsilon^{\rho}TK\nonumber\\
&  =\varepsilon+T(2\left\Vert \mu\right\Vert _{\mathcal{I}_{per}}+2^{\rho
}\varepsilon^{\rho-1}K)\varepsilon<2\varepsilon, \label{aux30}%
\end{align}
for all $u\in\mathcal{B}_{\varepsilon}$ and thus $\Phi(\mathcal{B}%
_{\varepsilon})\subset\mathcal{B}_{\varepsilon}.$ From Lemma \ref{lem-est1},
we\ also have%

\begin{align}
\left\Vert \Phi(u)-\Phi(v)\right\Vert _{L^{\infty}(-T,T;\mathcal{I}_{per})}
&  =\Vert L_{\mu,per}(u)-L_{\mu,per}(v)\Vert_{L^{\infty}(-T,T;\mathcal{I}_{per})}\nonumber\\
&+\left\Vert B_{per}(u)-B_{per}(v)\right\Vert _{L^{\infty}(-T,T;\mathcal{I}_{per})}\nonumber\\
&  \leq T\left\Vert \mu\right\Vert _{\mathcal{I}_{per}}\Vert u-v\Vert
_{L^{\infty}(-T,T;\mathcal{I}_{per})}\nonumber\\
&  +KT\Vert u-v\Vert_{L^{\infty}(-T,T;\mathcal{I}_{per})}\left(  \Vert
u\Vert_{L^{\infty}(-T,T;\mathcal{I}_{per})}^{\rho-1}+\Vert v\Vert_{L^{\infty
}(-T,T;\mathcal{I}_{per})}^{\rho-1}\right) \nonumber\\
&  \leq T\left(  \left\Vert \mu\right\Vert _{\mathcal{I}_{per}}+K2^{\rho
}\varepsilon^{\rho-1}\right)  \Vert u-v\Vert_{L^{\infty}(-T,T;\mathcal{I}_{per})}, \label{aux3}%
\end{align}
for all $u,v\in\mathcal{B}_{\varepsilon}.$ In view of (\ref{cond-small}),
(\ref{aux30}) and (\ref{aux3}), the map $\Phi$ is a contraction in
$\mathcal{B}_{\varepsilon}$ and then, the Banach fixed point theorem assures
the existence of a unique solution $u\in L^{\infty}(-T,T;\mathcal{I}_{per})$ for
(\ref{int2}) such that $\Vert u\Vert_{L^{\infty}(-T,T;\mathcal{I}_{per})}\leq2\Vert u_0\Vert_{\mathcal{I}_{per}}$.

On the other hand if $u,v$ are two solutions with respective initial data
$u_{0},v_{0}$ then, similarly in deriving (\ref{aux3}), one obtains
\begin{align*}
\Vert u-v\Vert_{L^{\infty}(-T,T;\mathcal{I}_{per})} &  \leq\Vert u_{0}-v_{0}\Vert_{L^{\infty
}(-T,T;\mathcal{I}_{per})}+\Vert L_{\mu,per}(u)-L_{\mu,per}(v)\Vert_{L^{\infty}(-T,T;\mathcal{I}_{per})}\nonumber\\
&+\left\Vert B_{per}(u)-B_{per}(v)\right\Vert _{L^{\infty
}(-T,T;\mathcal{I}_{per})}\\
&  \leq\Vert u_{0}-v_{0}\Vert_{L^{\infty}(-T,T;\mathcal{I}_{per})}+T\left(
\left\Vert \mu\right\Vert _{\mathcal{I}_{per}}+2^{\rho}\varepsilon^{\rho
-1}K\right)  \Vert u-v\Vert_{L^{\infty}(-T,T;\mathcal{I}_{per})},
\end{align*}
which, in view of (\ref{cond-small}), implies the desired local Lipschitz continuity.

\textbf{Proof of (2). }It follows by proceeding entirely parallel to the
proof of item (1) by replacing $\mathcal{I}_{per}$ by $\mathcal{I}.$ \fin

\section{Appendix}

Next we present a different proof of Theorem \ref{expli}. For $\sigma>0$, from Theorem 3.1 in Albeverio {\it et al.} \cite{AGHH} , the fundamental solution $S_\sigma(x, y;t)$ to the Schr\"odinger equation (\ref{SCH1}) is given by
\begin{equation}\label{kern}
S_\sigma(x, y,t)=S(x-y;t)-\frac{\sigma}{2}\int_0^\infty e^{-\frac{\sigma}{2} u} S(u+|x|+|y|;t)du
\end{equation}
where $S(x;t)$ is the free propagator in $\mathbb R$, i.e.
$$
S(x,t)= \frac{e^{-x^2/{4it}}}{(4i\pi t)^{1/2}},\quad t>0
$$
 and $e^{it\Delta}f(x)=S(x;t)\ast_x f(x)$. Then we have the representation
 \begin{equation}\label{propa}
 e^{-itH_\sigma}f(x)=\int_{\mathbb R} S_\sigma(x, y,t)f(y)dy.
  \end{equation}
Next, we consider $x>0$ and $f\in L^1(\mathbb R)$ with $supp f\subset (-\infty, 0]$. Then, since
$$
S(u+x-y;t)= (e^{-it\xi^2})^{\vee}(u+x-y)=(e^{-it\xi^2}e^{i\xi(x+u)})^{\vee}(y)
$$
we obtain via Parseval identity that  (\ref{propa})  can be re-write for $x>0$ in the form
 \begin{equation}\label{propa1}
 \begin{aligned}
 e^{-itH_\sigma}f(x)&=e^{it\Delta}f(x)\chi^0_+(x) +\int_{\mathbb R}e^{-ity^2}\Big[\int_{\mathbb R} -\frac{\sigma}{2} e^{\frac{\sigma}{2} s}\chi^0_-(s)e^{-iys}ds\Big] \widehat{f}(y)e^{iyx}dy\\
 &=e^{it\Delta}f(x)\chi^0_+(x) +\int_{\mathbb R}e^{-ity^2}\widehat{\rho_\sigma\ast f}(y)e^{iyx}dy\\
 &=e^{it\Delta}f(x)\chi^0_+(x) +e^{it\Delta}(\rho_\sigma\ast f)(x)\chi^0_+(x)
 \end{aligned}
 \end{equation}
where $\rho_\sigma(x)=-\frac{\sigma}{2} e^{\frac{\sigma}{2} s}\chi^0_-(s)$. Similarly, for $x<0$ we obtain the representation
 \begin{equation}\label{propa2}
e^{itH_\sigma}f(x)=e^{it\Delta}f(x)\chi^0_-(x) +e^{it\Delta}(\rho_\sigma\ast f)(-x)\chi^0_-(x).
\end{equation}
From (\ref{propa1})-(\ref{propa2}) we obtain the formula (\ref{pospro}).

For $\sigma<0$, Theorem 3.1 in  \cite{AGHH} establishes that the fundamental solution $S_\sigma(x, y;t)$ to the Schr\"odinger equation (\ref{SCH1}) is given by
$$
S_\sigma(x, y,t)=S(x-y;t) +e^{i\frac{\sigma^2 t}{4}}\Psi_\sigma(x)\Psi_\sigma(y)+\frac{\sigma}{2}\int_0^\infty e^{\frac{\sigma}{2} u} S(u-|x|-|y|;t)du.
$$
where $\Psi_\sigma$ is the (normalized) eigenfunction defined in Theorem \ref{resol5b}. Then, a similar analysis as above produces the formula (\ref{negpro}).

We note that since $|S(x;t)|\leqq C_0 t^{-1/2}$ for every $x\in \mathbb R$ and $t>0$ we obtain from (\ref{kern})  that $|S_\sigma(x,y;t)|\leqq 2C_0 t^{-1/2}$ for all $x, y\in \mathbb R$ and $t>0$. Therefore from (\ref{propa}) we obtain the dispersive estimate
$$
\|e^{-itH_\sigma}f\|_{\infty}\leqq 2C_0 t^{-1/2} \|f\|_1,
$$
which implies the estimate (\ref{stric2}) for $\sigma >0$.

\vskip0.2in

\textbf{ACKNOWLEDGEMENTS:} J. Angulo was partially supported by Grant CNPq/Brazil. L.C.F. Ferreira was supported by FAPESP-SP and CNPq/Brazil.

\end{document}